\documentclass{iopart}
\usepackage{iopams}
\usepackage{cite}
\usepackage{amssymb,amsfonts,euscript,graphicx,bm}
  \usepackage{epstopdf}
 \usepackage{graphics,setstack} 
 \usepackage[colorlinks=true]{hyperref}
 \hypersetup{urlcolor=blue, citecolor=red}
 \usepackage[latin1]{inputenc}

\eqnobysec

\newcommand{\z}{\mathbf{z}}

\newcommand{\W}{\mathbf{W}}
\newcommand{\p}{\mathbf{p}}
\newcommand{\q}{\mathbf{q}}

\newcommand{\Markov}[2]{\underset{#1}{\overset{#2}{\rightleftharpoons}}}

\renewcommand{\e}{{\rm e}}

\def\P{{\mathbb P}} 
 
\def\R{{\mathbb R}}
\newcommand{\LL}{{\cal L}}

\newcommand{\meas}{\mathcal{M}}

\newcommand{\mI}{\mathcal{J}}
\newcommand{\mF}{\mathbf{F}} %bold F
\newcommand{\mmu}{{\bm \mu}} %bold mu
\newcommand{\mA}{\mathbf{A}} %bold A
\newcommand{\mbp}{\mathbf{p}} % bold p
\newcommand{\mR}{{\mathbf R}} %bold R
\renewcommand{\theequation}{\arabic{section}.\arabic{equation}}

\begin{document}

\title{On the Hamiltonian structure of large deviations in stochastic hybrid systems}
\author{Paul C. Bressloff$^{1}$ and Olivier Faugeras$^2$}
\address{$^1$Department of Mathematics, University of Utah, 155 South 1400 East, Salt Lake City, Utah 84112, USA\\
 $^2$NeuroMathComp Group, INRIA, Sophia Antipolis, France}
\ead{bressloff@math.utah.edu,olivier.faugeras@inria.fr}
\date{\today}

\begin{abstract}

We develop the connection between large deviation theory and more applied approaches to stochastic hybrid systems by highlighting a common underlying Hamiltonian structure. A stochastic hybrid system involves the coupling between a piecewise deterministic dynamical system in $\R^d$ and a time-homogeneous Markov chain on some discrete space $\Gamma$. We assume that the Markov chain on $\Gamma$ is ergodic, and that the discrete dynamics is much faster than the piecewise deterministic dynamics (separation of time-scales). Using the Perron-Frobenius theorem and the calculus-of-variations, we evaluate the rate function of a large deviation principle in terms of a classical action, whose Hamiltonian is given by the Perron eigenvalue of a $|\Gamma|$-dimensional linear equation. The corresponding linear operator depends on the transition rates of the Markov chain and the nonlinear functions of the piecewise deterministic system. The resulting Hamiltonian is identical to one derived using path-integrals and WKB methods. We illustrate the theory by considering the example of stochastic ion channels. Finally, we indicate how the analysis can be extended to a multi-scale stochastic process, in which the slow dynamics is described by a piecewise SDE.
 \end{abstract}
 
 \maketitle

\section{Introduction}\label{sec:Intro}

There are a growing number of problems in biology that involve the coupling between a piecewise deterministic dynamical system in $\R^d$ and a time-homogeneous Markov chain on some discrete space $\Gamma$, resulting in a stochastic hybrid system, also known as a piecewise deterministic Markov process (PDMP) \cite{davis:84}. One important example is given by membrane voltage fluctuations in neurons due to the stochastic opening and closing of ion channels \cite{Fox94,Chow96,Keener11,Goldwyn11,Buckwar11,Pakdaman10,Wainrib12,NBK13,Bressloff14a}. Here the discrete states of the ion channels evolve according to a continuous-time Markov process with voltage-dependent transition rates and, in-between discrete jumps in the ion channel states, the membrane voltage evolves according to a deterministic equation that depends on the current state of the ion channels. In the limit that the number of ion channels goes to infinity, one can apply the law of large numbers and recover classical Hodgkin-Huxley type equations. However, finite-size effects can result in the noise-induced spontaneous firing of a neuron due to channel fluctuations. Other examples of stochastic hybrid systems include genetic switches \cite{Kepler01,Newby12}, motor-driven intracellular transport \cite{Newby10b,Newby11}, and stochastic neural networks \cite{Bressloff13}. 

In many of the above examples, one finds that the transition rates between the discrete states $n\in \Gamma$ are much faster than the relaxation rates of the piecewise deterministic dynamics for $x\in \R^d$. Thus there is a separation of time scales between the discrete and continuous processes, so that if $t$ is the characteristic time-scale of the Markov chain then $\epsilon t$ is the characteristic time-scale of the relaxation dynamics for some small positive parameter $\epsilon$. Assuming that the Markov chain is ergodic, in the limit $\epsilon\rightarrow 0$ one obtains a deterministic dynamical system in which one averages the piecewise dynamics with respect to the corresponding unique stationary measure. This then raises the important problem of characterizing how the law of the underlying stochastic process approaches this deterministic limit in the case of weak noise, $0<\epsilon \ll 1$.

A rigorous mathematical approach to addressing the above issue is {\em large deviation theory}, which has been developed extensively within the context of stochastic differential equations (SDEs) \cite{Freidlin98,Dembo98,Touchette09}. In particular, consider some random dynamical system in $\R^d$ for which there exists a well defined probability density functional or law $P_{\epsilon}[x]$ over the different sample trajectories $\{x(t)\}_0^T$ in a given time interval $[0,T]$. Here $\epsilon$ is a small parameter that characterizes the noise level, with $x(t)$ given by the solution $x^*(t)$ of some ODE $\dot{x}=F(x)$ in the limit $\epsilon \rightarrow 0$. A large deviation principle (LDP) for the random paths of the SDE over some time interval $[0,T]$ is 
\[P_{\epsilon}[x]\sim \e^{-J_T[x]/\epsilon},\quad \epsilon \rightarrow 0,\]
where $J_T[x]$ is known as the rate function and $J_T[x^*]=0$. In the case of SDEs, the rate function can be interpreted as a classical action with corresponding Lagrangian $L$ \cite{Freidlin98},
\[J_T[x]=\int_0^TL(x,\dot{x})dt.\] 
Such a Lagrangian formulation is more amenable to explicit calculations.
In particular,  it can be used to solve
various first passage time problems associated with the escape from a fixed point attractor of the underlying deterministic system in the weak noise limit. This involves finding the most probable paths of escape, which minimize the action with respect to the set of all trajectories emanating from the fixed point. Evaluating the action along a most probable path from the fixed point to another point $x$ generates a corresponding {\em quasipotential} $\Phi(x)$. From classical variational analysis, it can be shown that the quasipotential satisfies a Hamilton-Jacobi equation $H(x,\partial_x\Phi)=0$, where $H$ is the Hamiltonian obtained from $L$ \cite{Freidlin98} via a Fenchel-Legendre transformation:
\[H(x,p)=\underset{y}\sup\{py -L(x,y)\}.\]
The optimal paths of escape correspond to solutions of Hamilton's equations on the zero energy surface ($H=0$). Interestingly, the same quasipotential is obtained by considering a Wentzel-Kramers-Brillouin (WKB) approximation of the stationary state of the continuous process $x(t)$. More recently, large deviation theory has been applied to stochastic hybrid systems \cite{fagg08,fagg09,Kifer09}.

Independently of the developments in large deviation theory, a variety of techniques in applied mathematics and mathematical physics have been used to solve first passage time problems in biological applications of stochastic hybrid systems. These include WKB approximations and matched asymptotics \cite{Keener11,Newby11,Newby12,NBK13,Bressloff14a}, and path-integrals \cite{Bressloff14}. Although such approaches are less rigorous than large deviation theory, they are more amenable to explicit calculations. In particular, they allow one to calculate the prefactor in Arrhenius-like expressions for mean first passage times, rather than just the leading order exponential behavior governed by the quasipotential. The main aim of this paper is to make explicit the connection between large deviation theory and more applied approaches to stochastic hybrid systems, by highlighting the common underlying Hamiltonian structure. 

In \S 2 we define a stochastic hybrid system and specify our various mathematical assumptions. We then construct a WKB approximation of solutions to the differential Chapman-Kolmogorov equation, which describes the evolution of the probability density function for a stochastic hybrid system (\S 3). In particular, we show that the resulting quasipotential satisfies a Hamilton-Jacobi equation, whose associated Hamiltonian is given by the Perron eigenvalue of a corresponding matrix equation that depends on the transition rates of the Markov chain and the nonlinear functions of the piecewise deterministic system. This Hamiltonian is identical to the one rigorously derived by Kifer using large deviation theory \cite{Kifer09} and by Bressloff and Newby using formal path-integral methods \cite{Bressloff14,Bressloff15}. Further insights into the Hamiltonian structure of large deviations in stochastic hybrid systems are developed in \S 4 (and appendix A), where we present an alternative derivation of the Hamiltonian that does not require the full machinery of large deviation theory. We take as our starting point an LDP for stochastic hybrid systems due to Faggionato et al. \cite{fagg08,fagg09}, which is equivalent to the LDP of Kifer \cite{Kifer09}. (For completeness, we summarize a simplified version of Kifer's formulation in appendix B.) Using the Perron-Frobenius theorem and the calculus-of-variations, we evaluate the LDP rate function in terms of a classical action, whose equations of motion are given by the expected Hamiltonian dynamical system. We then consider some specific examples in order to illustrate our analysis (\S 5). Finally, in \S 6, we briefly indicate how the analysis can be extended to a multi-scale stochastic process, in which the slow dynamics is given by a piecewise SDE.

\section{The model}
Consider a one-dimensional stochastic hybrid system also called a piecewise deterministic Markov process (PDMP) \cite{davis:84,Kifer09}. Its states are described by a pair
$(x,n) \in \Omega \times \{1,\cdots,K\}$, where $x$ is a continuous variable in a connected, bounded domain $\Omega\subset \R^d$ with a regular boundary $\partial\Omega$   and $n$ a discrete internal state variable taking values in $\Gamma \equiv \{1,\cdots,K\}$. 
%(Note that one could easily extend our analysis to higher-dimensions, $x\in \R^d$. In this case $\Omega$ is taken to be a connected, bounded domain with a regular boundary $\partial\Omega$. However, for notational simplicity, we restrict ourselves to the case $d=1$.) 
When the internal state is $n$, the system evolves according to the ordinary differential equation (ODE)
\begin{equation}
\label{eq:deterministic}
\dot{x}=%\frac{1}{\tau} 
\mF_n(x),
\end{equation}
where the vector  field $\mF_n: \Omega \to \R^d$ is a continuous function, locally Lipschitz. That is, given a compact subset $\mathcal{K}$ of $\Omega$, there exists a positive constant $K_n(\mathcal{K})$ such that
\begin{equation}
|\mF_n(x)-\mF_n(y)|\leq K_n(\mathcal{K}) |x-y|,\quad \forall x,y\in \Omega.
\end{equation}
%Here $\tau$ is a fixed positive time constant that characterizes the relaxation rate of the $x$-dynamics. 
We assume that the dynamics of $x$ is confined to the domain $\Omega$ so that we have existence and uniqueness of a trajectory for each $n$. One final mild constraint is that the vector field does not have identical components anywhere in $\Omega$, that is, for any $x\in \Omega$, there exists at least one pair $(n,m)\in \Gamma, n\neq m$ for which $\mF_m(x)\neq \mF_n(x)$.

In order to specify how the system jumps from one internal state to the other for each $n \in \Gamma$, %we consider the positive time constant $\tau_n$ and 
the function $W_{nm}(x)$ defined on $\Gamma \times \Gamma \times \R$ with the following properties
\begin{enumerate}
\item $\forall n,\,m,\ W_{nm}(\cdot)$ is $C^1(\Omega)$
\item $\forall x,\,n,\ W_{n \cdot}(x)$ is a probability measure on $\Gamma$ such that $W_{nn}(x)=0$.
\end{enumerate}
At each point $x$ of $\Omega$ this defines a time-homogeneous Markov chain with state space $\Gamma$. We make the further assumption that this chain is irreducible (i.e. that the $K \times K$ matrix $\W(x)$ is irreducible) for all $x$ in $\Omega$, meaning that there is a non-zero probability of transitioning, possibly in more than one step, from any state to any other state of the Markov chain. This implies the existence of a unique invariant probability measure on $\Gamma$, noted $\rho(x,n)$, such that
\[
\sum_{n \in \Gamma} \rho(x,n) W_{nm}(x)=\rho(x,m)\ \forall m \in \Gamma,\,\forall x \in \Omega
\]
or in matrix form:
\begin{equation}\label{eq:invariant}
\rho(x,\cdot) \W(x) = \rho(x,\cdot),
\end{equation}
where $\rho(x,\cdot)$ is the $K$-dimensional row vector of coordinates $\rho(x,n)$, $n \in \Gamma$. The existence of the unique invariant measure can be seen as a consequence of the well-known Perron-Frobenius theorem\footnote{A finite-dimensional, real irreducible square matrix with non-negative entries has a unique largest positive real eigenvalue (the Perron eigenvalue) and the corresponding eigenvector has strictly positive components \cite{Grimmett}. 
 %The result also holds for a matrix with non-negative entries, provided that the matrix is irreducible. However, there can now be complex eigenvalues with the same absolute value as the Perron eigenvalue.
 }.

The hybrid evolution of the system is described as follows. Suppose the system starts at time zero in the state $(x_0,\,n_0)$. Call $x_0(t)$ the solution of (\ref{eq:deterministic}) with $n=n_0$ such that $x_0(0)=x_0$. Let $\theta_1$ be the random variable such that 
\[\P(\theta_1 > t) = \exp \left( -\frac{t}{\epsilon} \right).\]
Then in the random time interval $[0,\,\theta_1)$ the state of the system is $(x_0(s),n_0)$. We draw a value of $\theta_1$ from the corresponding probability density
\[p(t)=\frac{1}{\epsilon} \exp\left( -\frac{t}{\epsilon} \right).\]
We now choose an internal state $n_1 \in \Gamma$ with probability $W_{n_0n_1}(x_0(\theta_1))$ and call $x_1(t)$ the solution of the following Cauchy problem on $[\theta_1,\infty)$:
\[
\left\{
\begin{array}{lcl}
\dot{x}_1(t) & = & \mF_{n_1}(x_1(t)),\quad t \geq \theta_1\\ \\
x_1(\theta_1) & = & x_0(\theta_1)
\end{array}
\right.
\]
Iterating this procedure, we construct a sequence of increasing jumping times $(\theta_k)_{k \geq 0}$ (setting $\theta_0=0$) and a corresponding sequence of internal states $(n_k)_{k \geq 0}$. The evolution $(x(t),\,n(t))$ is then defined as
\begin{equation}
\label{hs}
(x(t),n(t))=(x_k(t),n_k) \quad \mbox{if}\ \theta_k \leq t < \theta_{k+1}.
\end{equation}
Note that the path $x(t)$ is continuous and piecewise $C^1$. In order to have a well-defined dynamics on $[0,T]$, it is necessary that almost surely the system makes a finite number of jumps in the time interval $[0,T]$. This is guaranteed in our case.

Given the above iterative definition of the stochastic hybrid process, let $X(t)$ and $N(t)$ denote the stochastic continuous and discrete variables, respectively, at time $t$, $t>0$, given the initial conditions $X(0)=x_0,N(0)=n_0$. Although the evolution of the continuous variable $X(t)$ or the discrete variable $N(t)$ is non-Markovian, it can be proven that the joint evolution $(X(t),N(t))$ is a strong Markov process \cite{davis:84}. The corresponding infinitesimal generator $L$ is easily shown to be equal to
\[
L g(x,n)=\mF_n(x) \cdot \nabla_x g(x,n)+\frac{1}{\epsilon} \sum_{m \in \Gamma} W_{nm}(x)(g(x,m)-g(x,n)),
\]
for regular functions $g: \Omega \times \Gamma \to \R$. We rewrite this as
\[
L g(x,n)=\mF_n(x) \cdot \nabla_x g(x,n)+\frac{1}{\epsilon} \sum_{m \in \Gamma} A_{nm}(x)g(x,m),
\]
where
\[
A_{nm}(x)=W_{nm}-\delta_{nm}.
\]
Note that $\sum_{m\in \Gamma}A_{nm}=0\quad  \forall n \in \Gamma$, and, in matrix form,
\begin{equation}\label{eq:AR}
\mA(x)=\W(x)-\mathbf{I}_K,
\end{equation}
where $\mathbf{I}_K$ is the $K \times K$ identity matrix. Note that (\ref{eq:invariant}) implies
\[
\rho(x,\cdot)\mA(x)=\mathbf{0}
\]

%If we let $\W^*$ be the vector of coordinates
%\[
%W^*_n=\sum_{k\in \Gamma}W_{kn}=1/\hat{\tau}_n,
%\]
%then we can write the matrix ${\bf A}$ as
%$
%\A=\W-{\rm diag}(\W^*)$.
%Although the evolution of the continuous variable $X(t)$ or the discrete variable $N(t)$ is non-Markovian, it can be proven that the joint evolution $(X(t),N(t))$ is a strong Markov process \cite{davis:84}. 

In this paper, we focus on the stochastic dynamics in the limit of fast kinetics, that is, $\epsilon \rightarrow 0$.

Let us now define the averaged vector field $\overline{\mF}: \R \to \R$ by
\[
\overline{\mF}(x)=\sum_{n \in \Gamma} \rho(x,n) \mF_n(x)
\]
It can be shown \cite{fagg08} that, given the assumptions on the matrix $\W$, the functions $\rho(x,n)$ on $\Omega$ belong to $C^1(\Omega)$ for all $n \in \Gamma$ and that this implies that $\overline{\mF}(x)$ is locally Lipschitz. Hence, for all $t \in [0,T]$, the Cauchy problem
\begin{equation}
\label{mft}
\left\{
\begin{array}{lcl}
\dot{x}(t) & = &  \overline{\mF}(x(t))\\
x(0) & = & x_0
\end{array}
\right.
\end{equation}
has a unique solution for all $x_0 \in \Omega$. Intuitively speaking, one would expect the stochastic hybrid system (\ref{eq:deterministic}) to reduce to the deterministic dynamical system (\ref{mft}) in the limit $\epsilon \rightarrow 0$. That is, for sufficiently small $\epsilon$, the Markov chain undergoes many jumps over a small time interval $\Delta t$ during which $\Delta x\approx 0$, and thus the relative frequency of each discrete state $n$ is approximately $\rho(x,n)$. This can be made precise in terms of a Law of Large Numbers for stochastic hybrid systems proven in \cite{fagg08}.

\section{WKB method and quasipotentials}\label{sect:WKB}
For notational simplicity, we restrict ourselves in this section to the case $d=1$.

Introduce the probability density $\rho(x,n,t|x_0,n_0,0) $ with
\[\P\{X(t)\in (x,x+dx),\, N(t) =n|x_0,n_0)=\rho(x,n,t|x_0,n_0,0)dx.\]
%We also fix the units of time by setting $\tau=1$ and introducing the scalings
%\begin{equation}
%\tau_n=\epsilon \hat{\tau}_n,\quad W_{mn}\rightarrow W_{mn}/\hat{\tau}_n,
%\end{equation}
%with $\hat{\tau}_n$ independent of $\epsilon$. 
It follows that $\rho$ evolves according to the forward differential Chapman-Kolmogorov (CK) equation \cite{Gardiner09,Bressloff14}
\begin{equation}
\label{CK}
\frac{\partial \rho}{\partial t}=L_t\rho ,
\end{equation}
with (dropping the explicit dependence on initial conditions)
\begin{equation}
L_t\rho(x,n,t)=-\frac{\partial F_n(x)\rho(x,n,t)}{\partial x}+\frac{1}{\epsilon}\sum_{m\in \Gamma}A^{\top}_{nm}(x)\rho(x,m,t),
\end{equation}

\begin{figure}[t!]
\begin{center}
\includegraphics[width=10cm]{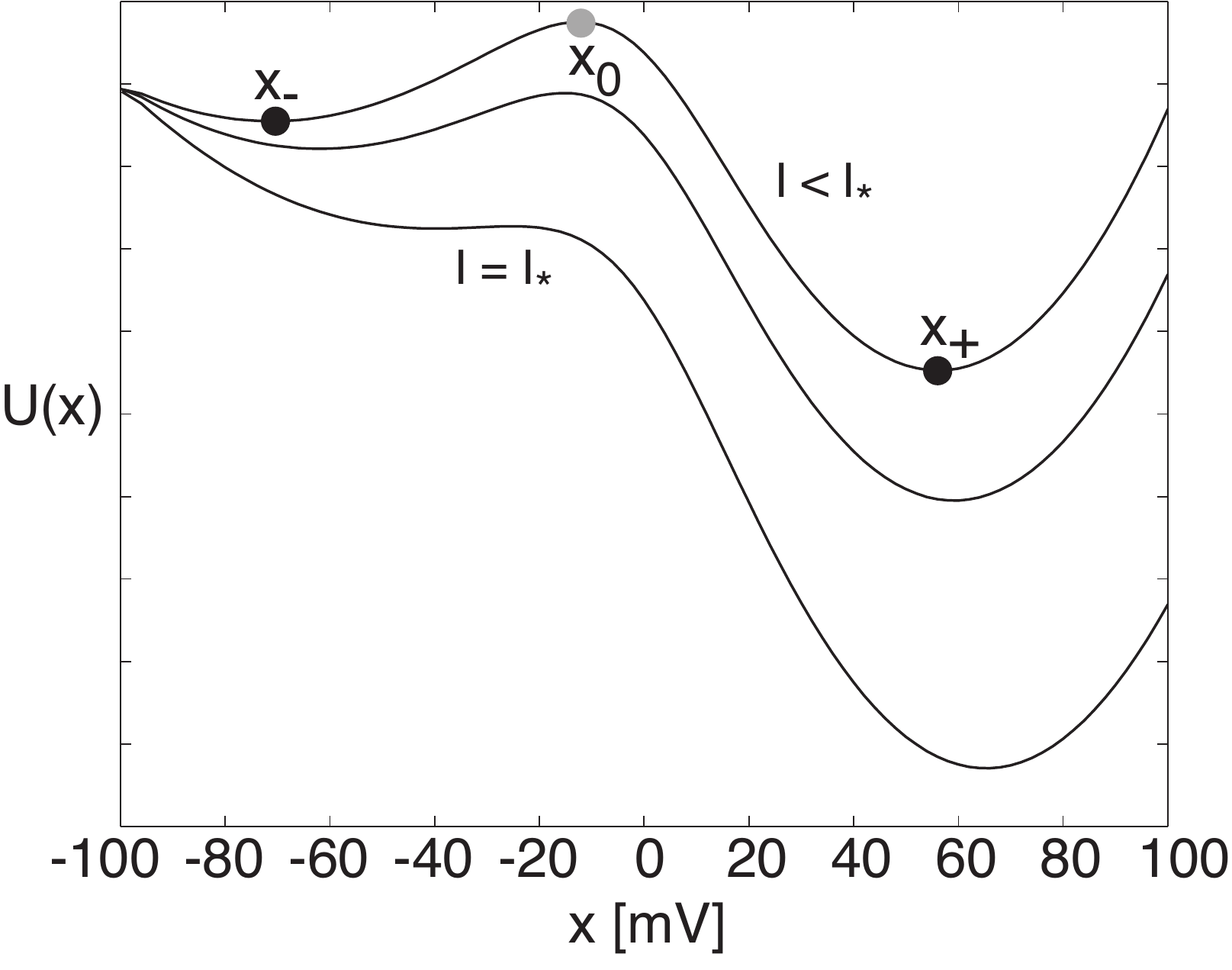}
\caption{\small Plot of deterministic potential $U(x)$ as a function of voltage $x$ for different values of the external stimulus current $I$. Parameter values are $V_{\rm Na}=120$ mV, $V_{L}=-62.3$ mV, $g_{\rm Na}=4.4$ mS/cm$^2$, $g_{L}= 2.2$ mS/cm$^2$, and $\alpha(x)=\beta \exp[(x-v_1)/v_2]$ with $\beta=0.8$ s$^{-1}$, $v_1=-1.2mV$, $v_2=18mv$.}
\label{wellx}
\end{center}
\end{figure}

Suppose that the averaged equation (\ref{mft}) is bistable with a pair of stable fixed points $x_{\pm}$ separated by an unstable fixed point $x_0$ (see section 4 and Fig. \ref{wellx} for an explicit example). Assume that the stochastic system is initially at $x_-$. On short time scales ($t \ll 1/\epsilon$) the system rapidly converges to a quasistationary solution within the basin of attraction of $x_-$, which can be approximated by a Gaussian solution of the reduced FP equation obtained using a quasi-steady-state (QSS) diffusion or adiabatic approximation of the CK equation (\ref{CK}) \cite{Newby10a}. However, on longer time scales, the survival probability slowly decreases due to rare transitions across $x_0$ at exponentially small rates, which cannot be calculated accurately using the  diffusion approximation. One thus has to work with the full CK equation (\ref{CK}) supplemented by an absorbing boundary conditions at $x_0$:
\begin{equation}
\rho(x_0,n,t)=0, \, \mbox{for all } \, n \in \Sigma,
\end{equation}
where $\Sigma\subset \Gamma$ is the set of discrete states $n$ for which $F_n(x_0)<0$.
The initial condition is taken to be
\begin{equation}
\label{icon}
\rho(x,n,0)=\delta(x-x_-)\delta_{n,\bar{n}}.
\end{equation}
Let $T$ denote the (stochastic) first passage time for which the system first reaches $x_0$, given that it started at $x_-$. The distribution of first passage times is related to the survival probability that the system hasn't yet reached $x_0$, that is, 
\[
\mbox {Prob}\{t>T\}=S(t)\equiv \int_{-\infty}^{x_0} \sum_{n\in \Gamma} \rho(x,n,t)dx .
\]
The first passage time density is then
\begin{equation}
f(t)=-\frac{dS}{dt}=-\int_{-\infty}^{x_0}\sum_{n \in \Gamma} \frac{\partial \rho(x,n,t)}{\partial t}dx .
\end{equation}
Substituting for $\partial \rho /\partial t$ using the CK equation (\ref{CK}) shows that
\begin{eqnarray}
\fl f(t)&=&\int_{-\infty}^{x_0} \left [\sum_{n \in \Gamma}\frac{\partial [F_n(x) \rho(x,n,t)]}{\partial x}\right ]dx= \sum_{n \in \Gamma}F_n(x_0)\rho(x_0,n,t) .
\label{fTP}
\end{eqnarray}
We have used $\sum_{n}{A}^{\top}_{nm}(x)=0$ and $\lim_{x\rightarrow-\infty}F_n(x)\rho(x,n,t)=0$. The first passage time density can thus be interpreted as the probability flux $J(x,t)$ at the absorbing boundary, since we have the conservation law
\begin{equation}
\sum_{n \in \Gamma}\frac{\partial \rho(x,n,t)}{\partial t}=-\frac{\partial J(x,t)}{\partial x},\quad J(x,t)=\sum_{n \in \Gamma}F_n(x)\rho(x,n,t).
\end{equation}

Suppose that the solution to the CK equation (\ref{CK}) with an absorbing boundary at $x=x_0$ has the eigenfunction expansion
\begin{equation}
\rho(x,n,t)=\sum_{r\geq 0}C_r\e^{-\mu_rt}V_r(x,n),
\end{equation}
where $\mbox{Re}[\mu_r]>0$ for all $r\geq 0$ and $(\mu_r,V_r)$ are an eigenpair of the linear operator ${\mathcal L}$ on the right-hand side of the CK equation (\ref{CK}):
\begin{equation}
\label{Lhat}
\fl {\mathcal L}V_r(x,n) \equiv \frac{d}{dx}(F_n(x)V_r(x,n))-\frac{1}{\epsilon}\sum_{m\in \Gamma}A^{\top}_{nm}V_r(x,m)=\mu_rV_r(x,n),
\end{equation}
together with the boundary conditions
\begin{equation}
V_r(x_0,n)=0,\, \mbox{for}\, n\in \Sigma .
\end{equation}
Furthermore, assume that there exists a simple, real eigenvalue eigenvalue $\mu_0$ such that $0<\mu_0<\mbox{Re}[\mu_1]\leq \mbox{Re}[\mu_2]\leq \ldots$ with $\mu_0\sim \e^{-C/\epsilon}$, whereas $\mbox{Re}[\mu_r]=O(1)$ for $r>0$. It follows that at large times we have the quasistationary approximation
\index{quasistationary approximation}
\begin{equation}
\label{qa}
\rho(x,n,t)\sim C_0 \e^{-\mu_0 t}V_0(x,n) .
\end{equation}
Substituting such an approximation into equation (\ref{fTP}) implies that $f(t)\sim \mu_0\e^{-\mu_0 t}$ with $\mu_0^{-1}$ the mean first passage time (MFPT). As shown elsewhere \cite{Keener11,Newby11}, one can determine $\mu_0$ by first using a WKB approximation to construct a quasistationary solution $V^{\epsilon}(x,n)$ for which ${\mathcal L}V^{\epsilon}=0$ and $V^{\epsilon}(x_0,n)\sim O(\e^{-\tau/\epsilon})$, and then performing an asymptotic expansion in order to match the quasistationary solution with the solution in a neighborhood of $x_0$. 
The WKB approximation $V^{\epsilon}(x)$ takes the form
\index{quasi-potential}
\begin{equation}
\label{WKBsh}
V^{\epsilon}(x,n)\sim \eta(x,n)\exp\left (-\frac{\Phi(x)}{\epsilon}\right ),
\end{equation}
where $\Phi(x)$ is the quasipotential. Substituting into the time-independent version of equation (\ref{CK}) yields
\begin{eqnarray}
\label{AAA}
&&\sum_{m \in \Gamma}\left (A^{\top}_{nm}(x)+\Phi'(x)\delta_{nm}F_m(x)\right )\eta(x,m) =\epsilon \frac{dF_n(x)\eta(x,n)}{d x},
\end{eqnarray}
where $\Phi'=d\Phi/dx$.
Introducing the asymptotic expansions $\eta \sim \eta^{(0)}+\epsilon \eta^{(1)}+\ldots$ and $ \Phi\sim \Phi_0+\epsilon \Phi_1+\ldots$, the leading order equation is
\begin{equation}
\sum_{m \in \Gamma} A^{\top}_{nm}(x)\eta^{(0)}(x,m)=- \Phi_0'(x)F_n(x)\eta^{(0)}(x,n).
\label{R0}
\end{equation}

It is well known that for SDEs in the weak-noise limit, the associated quasipotential satisfies a Hamilton-Jacobi equation \cite{Freidlin98,Stein97}. An analogous result also holds for jump Markov processes \cite{Hanggi84,Matkowsky85,Dykman94}. We now show that equation (\ref{R0}) can be reformulated as a Hamilton-Jacobi equation for the quasipotential $\Phi_0$ of a stochastic hybrid system.
First, we introduce the family of eigenvalue equations
\begin{equation}
\label{Perm0}
\fl \sum_{m \in \Gamma} A^{\top}_{nm}(x)z_m(x,p)+pF_n(x)z_n(x,p)=\lambda(x,p) z_n(x,p),
\end{equation}
which are parameterized by the pair $(x,p)$ with $p$ an auxiliary variable. An application of the Perron-Frobenius theorem along similar lines to  \ref{app:A} establishes that, for each fixed $(x,p)$, there exists a unique (Perron) eigenvalue whose associated eigenvector is strictly positive. Therefore, in the following we take $\lambda(x,p)$ to be the Perron eigenvalue and $z_n(x,p), n=1,\ldots,|\Gamma|,$ to be the corresponding eigenvector. Comparison of equation (\ref{R0}) with (\ref{Perm0}) then shows that we can make the identifications $\Phi_0'(x)= p$, $\eta^{(0)}(x,m)=z_m(x,p)$ and $\lambda(x,p)= 0$. That is, the quasipotential is the solution of the Hamilton-Jacobi equation
\begin{equation}
\label{HJ}
\lambda(x,\Phi_0'(x))=0.
\end{equation}
This is equivalent to finding zero energy solutions of Hamilton's equations 
\[
  \dot{x} = \frac{\partial \lambda(x,p)}{\partial p}, \quad  \dot{p} =- \frac{\partial \lambda(x,p)}{\partial x},
\]
and identifying $\Phi_0$ as the action along the resulting solution curve $(x(t),p(t))$:
\[
\Phi_0(x)=\int_{\bar{x}}^{x} p(x') dx'.
\]

Given the quasi-potential $\Phi_0$, the mean first passage time $\tau$ to escape from the fixed point at $x_-$ can be calculated by considering higher order terms in the WKB approximation, and using matched asymptotics to deal with the absorbing boundary at $x_0$. One finds that $\tau$ takes the general Arrhenius form \cite{Keener11,Newby11}
\[
\tau\sim \frac{\chi(x_0,x_-)}{\sqrt{|\Phi_0^{''}(x_0)|\Phi_0^{''}(x_-)}}\e^{[\Phi_0(x_0)-\Phi_0(x_-)]/\epsilon},
\]
where $\chi$ is an appropriate prefactor. Hence, the WKB method provides a powerful calculational tool. On the other hand, there is no {\em a priori} justification for interpreting the quasipotential and its associated Hamiltonian in terms of an underlying variational problem for optimal paths in the space of stochastic trajectories. This becomes crucial when solving escape problems in higher dimensions, since a metastable state is now surrounded by a non-trivial boundary (rather than a single point) and one needs to determine the relative weighting of optimal paths crossing different points on the boundary.  

\section{Classical action from a large deviation principle}

A powerful methodology for rigorously deriving a variational principle for optimal paths is large deviation theory \cite{Freidlin98}. The particular application of large deviation theory to stochastic hybrid systems has been developed in considerable detail by Kifer \cite{Kifer09}, who shows how to derive the action functional with Hamiltonian given by the Perron eigenvalue $\lambda(x,p)$ from an underlying large deviation principle (LDP). However, the technical nature of the analysis makes it difficult for a non-expert to understand. Therefore, in this section and appendix A we present an alternative derivation of the Hamiltonian system starting from the LDP of Faggionato et al. \cite{fagg08,fagg09}. The advantage of the latter LDP, which is equivalent to the LDP of Kifer (see \cite{Kifer09} and appendix B), is that it is relatively straightforward to write down. As far as we are aware, our derivation of the Hamiltonian action functional is new, and it provides additional insights into the origin and meaning of the ``momentum'' variable $p$. Finally, note that yet another method for deriving the action functional is to use formal path-integral methods \cite{Bressloff14}. One of the useful features of the path-integral formulation is that it can also be used to implement a variety of approximation schemes, including diffusion approximations and perturbation expansions of moment equations \cite{Bressloff15}.

In order to write down the LDP of Faggionato et al \cite{fagg08,fagg09}, it is first necessary to introduce some notation. Let $ \meas_+([0,T])$ denote the space of non-negative finite measures on the interval $[0,T]$ with a density w.r.t to the Lebesgue measure. It is convenient to isolate a special subset ${\mathcal Y}$ of $C([0,T],\Omega) \times \meas_+([0,T])^\Gamma$ as follows. Consider the pairs $\{(x(t),\psi(t))\}_{t\in [0,T]}$ such that the $K$-dimensional vector $\{\psi(t)\}_{t\in [0,T]}$ is an element of the product space $ \meas_+([0,T])^{\Gamma}$ satisfying the following extra condition  
\[\sum_{n\in\Gamma}\psi_n(t)=1,\]
where we have identified the measures $\psi_n(t)dt$ with their densities $\psi_n(t)$.
%In other words, for each $t\in [0,T]$, $\psi(t)=(\psi_1(t),\ldots,\psi_K(t))$ such that
%\psi_n(t)\geq 0,\quad 

A particular realization of the stochastic process, $\{(x(t),n(t))\}_{t\in [0,T]}$, then lies in ${\mathcal Y}$ with
\begin{equation}
\label{pip}
\psi_n(t)=1_{\{n(t)=n\}}\equiv \left \{ \begin{array}{cc} 1, & \mbox{ if } n(t)=n,\\ 0 , & \mbox{ if } n(t)\neq n \end{array} \right .
\end{equation}
and
\begin{equation}
\label{xx0}
x(t)=x_0+\int_0^T\sum_{n\in \Gamma} \psi_n(s)\mF_n(x(s))ds.
\end{equation}
 ${\mathcal Y}$  contains both the set of trajectories of the stochastic hybrid system with $\psi_n(t)$ given by equation (\ref{pip}) and $n(t)$ evolving according to the Markov chain, and the solution $x^*(t)$ of the averaged system (\ref{mft}) for which $\psi_n=\psi^*_n$ with $\psi_n^*(t)=\rho(x^*(t),n)$. It can be proven that ${\mathcal Y}$ is a compact subspace of $C([0,T],\Omega) \times \meas_+([0,T])^\Gamma$ with topology defined by the metric \cite{fagg08}
\begin{eqnarray}
 &d(\{(x(t),n(t))\}_{t\in [0,T]},\{(\widetilde{x}(t),\widetilde{n}(t))\}_{t\in [0,T]})\nonumber \\
&\qquad =\underset{t\in [0,T]}\sup|x(t)-\widetilde{x}(t)|+\sum_{n\in \Gamma}\underset{0\leq t\leq T }\sup \left |\int_0^t[\psi_n(s)-\widetilde{\psi}_n(s)]ds\right |.
\end{eqnarray}
Finally, we take $P^{\epsilon}_{x_0}$ to be the probability density functional or law of the set of trajectories $\{x(t)\}_{t\in [0,T]}\in C([0,T],\Omega)$. 

The following large deviation principle then holds \cite{fagg08,fagg09}:
Given an element $\{x(t)\}_{t\in [0,T]}\in C([0,T],\Omega)$,
\begin{eqnarray}
\label{LDP2}
\P^{\epsilon}_{x_0}\left [\{{x}(t)\}_{t\in [0,T]} \right ]\sim \e^{-J_T(\{x(t)\}_{t\in [0,T]})/\epsilon }.
\end{eqnarray}
where the rate function $J_{[0,T]}:  C([0,T],\Omega) \to [0,\infty)$ is given by
\begin{eqnarray}
\label{eq:rate2}
 &J_T(\{x(t)\}_{t\in [0,T]})\\
 &\qquad =\underset{\{\psi(t)\}_{t\in [0,T]}:\{x(t),\psi(t)\}_{t\in [0,T]} \in {\mathcal Y}_{x_0}}\inf  \int_0^T j(x(t),\psi(t))\,dt\nonumber
\end{eqnarray}
and 
\begin{equation}\label{eq:j}
j(x,\psi)=\sup_{\mR \in (0,\infty)^\Gamma} \sum_{(n,m) \in \Gamma\times \Gamma} \psi_n W_{nm}(x)\left[1-\frac{R_{m}}{R_n}\right].
\end{equation}
Here the symbol $\sim$ means asymptotic logarithmic equivalence in the limit $\epsilon \rightarrow 0$.

A key idea behind the above LDP is that a slow dynamical process coupled to the fast Markov chain on $\Gamma$ rapidly samples the different discrete states of $\Gamma$ according to some non-negative measure $\psi$. In the limit $\epsilon \rightarrow 0$, one has $\psi \rightarrow \rho$, where $\rho$ is the invariant measure of the Markov chain defined by (\ref{eq:invariant}). On the other hand, for small but non-zero $\epsilon$, $\psi$ is itself distributed according to a LDP, whereby one averages the different functions $\mF_n(x)$ over the measure $\psi$ to determine the dynamics of the slow system. In most biological applications, one is not interested in the internal discrete state of the system, that is, one only observes the statistical behavior of the continuous variable $x(t)$. For example, $x(t)$ could represent the membrane voltage of a neuron \cite{NBK13} or the position of a molecular motor along a microtubular track \cite{Newby10b,Newby11}. Faggionato et al. \cite{fagg08,fagg09} explicitly calculated the rate function (\ref{eq:rate2}) for a restricted class of stochastic hybrid systems, whose stationary density is exactly solvable. One major constraint on this class of model is that the vector field of the piecewise deterministic system has non-vanishing components within a given confinement domain. However, this constraint does not hold for biological systems such as ion channels  \cite{Keener11,NBK13,Bressloff14a}, gene networks \cite{Kepler01,Newby12}, and neural networks  \cite{Bressloff13} (Faggionato et al. were motivated by a model of molecular motors that is exactly solvable. Such a solvability condition also breaks down for molecular motors when local chemical signaling is taken into account \cite{Newby10b}.)

In appendix A we show (without the restrictions of Faggionato et al. \cite{fagg08,fagg09}) that the rate function $J_T(\{x(t)\}_{t\in [0,T]})$ of the non-Markov process $\{x(t)\}_{t\in [0,T]}$ can be written in the form 
of an action
\begin{equation}
\label{action}
J_T(\{x(t)\}_{t\in [0,T]})=\int_0^TL(x,\dot{x})dt,
\end{equation}
with Lagrangian given by
\[
L(x,\dot{x}) =\langle \mmu(x,\dot{x}),\,\dot{x}\rangle-\lambda(x(t),\mmu(x,\dot{x})),
\]
Here $\lambda(x,\mmu)$ is the Perron eigenvalue of the linear equation (\ref{Perm0}) (suitably generalized, see \ref{app:A}) for $p\rightarrow \mu$, that is,
\begin{equation}
\label{Perms}
%\sum_{m \in \Gamma} A_{nm}(x)R_m(x,\mu)+\mu F_n(x)R_n(x,\mu)=\lambda(x,\mu) R_n(x,\mu),
{\bf A}(x){\bf R}(x,\mmu)+({\bf F}(x){\bm \mu}) \circ{\bf R}(x,\mmu) =\lambda(x,\mmu){\bf R}(x,\mmu).
\end{equation}
Here, for any ${\bf a},{\bf b}\in \R^K$,
 \begin{equation}\label{eq:ab}
 [ {\bf a}\circ{\bf b}]_n \equiv [{\rm diag}({\bf a}){\bf b}]_n=a_nb_n,\,n=1,\cdots,K,
 \end{equation} 
the $K$-dimensional vector ${\mathbf F}(x){\bm \mu}$ is the product of the $K \times d$ matrix ${\mathbf F}(x)$ whose $K$ rows are the $d$-dimensional vectors $\mF_m(x)$, $m=1,\cdots,K$, with the $d$-dimensional vector $\bm \mu$, and
$\mmu=\mmu(x,\dot{x})$ is the solution of the invertible equation
\begin{equation}
\label{dotx0}
\dot{x}=\sum_{m \in \Gamma}\psi_m(x,\mmu)\mF_m(x),
\end{equation}
with
\begin{equation} \label{eq:psi}
 \psi_m(x,\mmu)= z_m(x,\mmu)R_m(x,\mmu), 
 \end{equation}
where $\z$ is the positive eigenvector of the adjoint equation
\begin{equation}
\label{Permz}
 \mA^\top(x) \z(x,\mmu)+ ({\bf F}(x){\bm \mu}) \circ \z(x,\mmu) = \lambda(x,\mmu)\z(x,\mmu)
 %\sum_{n \in \Gamma} z_n(x,\mu) A_{nm}(x)+\mu z_m(x,\mu)F_m(x)=\lambda(x,\mu) z_m(x,\mu)
\end{equation}
under the normalizations $\sum_n z_n=1$ and  $\langle \z,\,\mR \rangle=\sum_mz_m(x,\mu)R_m(x,\mu)=1$. 

Given the Lagrangian $L$, we can determine the Hamiltonian $H$ according to the Fenchel-Legendre transformation
\[
H(x,\mbp)=\underset{y \in \R^d}\sup \left [\langle \mbp-\mmu(x,y),\,y \rangle+\lambda(x,\mu(x,y))\right ].
\]
Evaluating the right-hand side yields the equation
\begin{eqnarray}
\label{var}
\mbp-\mmu(x,y)+ \frac{\partial \mmu}{\partial y}\left [\frac{\partial \lambda}{\partial \mmu}^\top-y\right ]=0
\end{eqnarray}
with
\[y=\sum_{m\in \Gamma}\psi_m(x,\mmu)\mF_m(x).\]
Differentiating equation (\ref{Perms}) with respect to $\mmu$ gives
\begin{equation}\label{eq:LR}
%\label{googl}
{\bf L}(x,{\bm \mu})\frac{\partial {\bf R}}{\partial {\bm \mu}}   ={\bf R}\frac{\partial \lambda}{\partial {\bm \mu}} -{\rm diag}(\mR) \mF(x),
\end{equation}
with ${\bf L}(x,{\bm \mu})=\left . {\bf L}(x,\p)\right |_{p_m=\langle {\bm \mu},\, \mF_m \rangle}$.
%\begin{eqnarray}
%\label{diff}
%& \sum_{m \in \Gamma} A_{nm}(x)\frac{\partial R_m(x,\mu)}{\partial \mu}+[\mu F_n(x)-\lambda(x,\mu))\frac{\partial R_n(x,\mu)}{\partial \mu}\nonumber \\
% & \qquad =\left [\frac{\partial \lambda(x,\mu)}{\partial \mu}-F_n(x)\right ] R_n(x,\mu).
%\end{eqnarray}
Since the adjoint of the linear operator on the left-hand side has a one-dimensional null space spanned by $\z$, see (\ref{Permz}), multiplying both sides of the equation by $\z^\top$, it follows from the normalization $\sum_mz_mR_m=1$ and (\ref{eq:psi}) that
\[\frac{\partial \lambda(x,\mmu)}{\partial \mmu}^\top=\sum_{m\in \Gamma}\psi_m(x,\mmu)\mF_m(x)=y.\] 
Equation (\ref{var}) thus shows that $\mbp=\mmu$, and we can identify $\mbp$ as the ``conjugate momentum'' of the Hamiltonian 
\begin{equation}
\label{H}
H=\lambda(x,\mbp),
\end{equation}
where $\lambda(x,\mbp)$ is the Perron eigenvalue of the linear equation (\ref{Perms}).

\section{Examples}

In this section we illustrate the Hamiltonian structure of stochastic hybrid systems by considering a few explicit models. 

\subsection{Binary model}
 We begin with the simple example of a binary stochastic hybrid process in dimension 1 ($d=1$) with two discrete states $n=0,1$, which was analyzed in some detail by Faggionato et al \cite{fagg09} using a different method. The latter authors exploited the fact that the model is exactly solvable, in the sense that the stationary density of the corresponding Chapman-Kolmogorov equation can be computed explicitly, and used a fluctuation-dissipation theorem to determine the Hamiltonian and quasipotential. Here we will obtain the same results more directly by calculating the Perron eigenvalue. One biological application of the binary model is to the bidirectional transport of a molecular motor along a one-dimensional microtubular track, in which $x$ represents the spatial location of the motor on the track and the two discrete states represent the motor moving either towards the $+$ end or $-$ end of the track. (In more complex models, the discrete space $\Gamma$ represents multiple internal conformational states of the motor, each of which has an associated velocity on the track \cite{Bressloff13}.)

Suppose that the continuous variable evolves according to piecewise dynamics on some finite interval $(a,b)$,
\begin{equation}
\dot{x}=F_n(x),\quad n=0,1,
\end{equation}
with $F_0,F_1$ continuous and locally Lipschitz. Suppose that $F_0,F_1$ are non-vanishing within the interval $(a,b)$, and $F_n(a)\geq  0,F_n(b) \leq  0$ for $n=0,1$; the dynamics is then confined to $(a,b)$. Denote the transition rates of the two-state Markov chain by $\omega_{\pm}(x)$ with
\[\{n=0\}\Markov{\omega_-(x)}{\omega_+(x)}\{n=1\}.\]
The stationary measure of the Markov chain is given by
\begin{equation}
\rho(x,0)=\frac{\omega_-(x)}{\omega_-(x)+\omega_-(x)},\quad \rho(x,1)=\frac{\omega_+(x)}{\omega_-(x)+\omega_-(x)}.
\end{equation}
The linear equation (\ref{Perms}) can be written as the two-dimensional system
\begin{equation}
\fl \left (\begin{array}{cc} -\omega_+(x)+pF_0(x)& \omega_+(x)\\ \omega_-(x) & -\omega_-(x)+pF_1(x) \end{array}\right )\left (\begin{array}{c} R_0 \\ R_1\end{array}\right )=\lambda \left (\begin{array}{c} R_0 \\ R_1\end{array}\right ).
\end{equation}
The corresponding characteristic equation is
\begin{eqnarray*}
0&=\lambda^2+\lambda[\omega_+(x)+\omega_-(x)-p(F_0(x)+F_1(x))]\\
&\quad +(pF_1(x)-\omega_-(x))(pF_0(x)-\omega_+(x))-\omega_-(x)\omega_+(x).
\end{eqnarray*}
It follows that the Perron eigenvalue is given by
\begin{eqnarray}
\lambda(x,p)=\frac{1}{2}\left [  \Sigma(x,p)+\sqrt{\Sigma(x,p)^2- 4\gamma(x,p)}    \right ]
\end{eqnarray}
where
\[\Sigma(x,p)=p(F_0(x)+F_1(x))-[\omega_+(x)+\omega_-(x)],\]
and
\[\gamma(x,p)=(pF_1(x)-\omega_-(x))(pF_0(x)-\omega_+(x))-\omega_-(x)\omega_+(x).\]
A little algebra shows that
\[D(x,p)\equiv \Sigma(x,p)^2- 4\gamma(x,p)=[p(F_0-F_1)-(\omega_+-\omega_-)]^2+\omega_+\omega_- >0\]
so that as expected $\lambda$ is real. Hence, from Hamilton's equations
\begin{eqnarray}
\fl \dot{x}&=\frac{\partial \lambda(x,p)}{\partial p} \nonumber \\
\fl &=\frac{F_0(x)+F_1(x)}{2}+\frac{\partial D(x,p)}{\partial p}\frac{1}{2\sqrt{D(x,p)}}\nonumber \\
\fl &=\frac{F_0(x)+F_1(x)}{2}+\frac{F_0(x)-F_1(x)}{2}\frac{p(F_0-F_1)-(\omega_+-\omega_-)}{\sqrt{[p(F_0-F_1)-(\omega_+-\omega_-)]^2+\omega_+\omega_- }}
\end{eqnarray}
which is the same result as obtained in example 10.4 of Faggionato et al \cite{fagg09}. Moreover,
writing
\[\dot{x}=F_0(x)\psi_0(x)+F_1(x)\psi_1(x),\]
we see that
\begin{equation}
\psi_0(x)=\frac{1}{2}\left [1+\frac{p(F_0-F_1)-(\omega_+-\omega_-)}{\sqrt{[p(F_0-F_1)-(\omega_+-\omega_-)]^2+\omega_+\omega_- }}
\right ],
\end{equation}
and
\begin{equation}
\psi_1(x)=\frac{1}{2}\left [1-\frac{p(F_0-F_1)-(\omega_+-\omega_-)}{\sqrt{[p(F_0-F_1)-(\omega_+-\omega_-)]^2+\omega_+\omega_- }}\right ],
\end{equation}
so that $\psi_{0,1}(x)\geq 0$ with $\psi_0(x)+\psi_1(x)=1$.

\subsection{Stochastic ion channels}

Another important example of stochastic hybrid systems is a conductance-based model of a neuron, in which the stochastic opening of membrane ion channels generates a stochastic ionic current that drives the membrane voltage. It is then possible that
ion channel noise induces spontaneous action
potentials (SAPs), which can have a large effect on a neuron's
function \cite{Fox94,Chow96,Keener11,Goldwyn11,Buckwar11,NBK13,Bressloff14a}. If SAPs are too frequent, a neuron cannot
reliably perform its computational role. Hence, ion channel
noise imposes a fundamental limit on the density of
neural tissue. Smaller neurons must function with fewer
ion channels, making ion channel fluctuations more significant
and more likely to cause a SAP. Here we will consider the simple case of a single type of ion channel, namely, a fast sodium (Na) channel, which was previously analyzed using WKB methods \cite{Keener11}. Let $x(t)$ denote the membrane voltage of the neuron at time $t$ and $N$ be the fixed number of sodium channels.
We assume that each channel can either be open $(O)$ or closed $(C)$, and can switch between each state according to the kinetic scheme
\begin{equation}
  \label{eq:4}
  C\Markov{\alpha(x)}{\beta(x)}O.
\end{equation}
with voltage-dependent transition rates.
(A more detailed biophysical model would need to treat each ion channel as a cluster of subunits rather than a single unit. In other words, the Markov chain of events associated with opening and closing of an ion channel would  involve transitions between more than two internal states,.)
The stochastic membrane voltage is taken evolves according to the piecewise deterministic equation
\begin{equation}
\label{pwv}
{\frac{dx}{dt}=F_n(x)\equiv \frac{n}{N}f(x)-g(x),}
\end{equation}
where 
\[f(x)=g_{\rm Na}(V_{\rm Na}-x), \quad g(x)=-g_{\rm L}[V_{L}-x]-I. \]
Here $g_{\rm Na}$ is the maximal conductance of a sodium channel and $V_{\rm Na}$ is the corresponding membrane reversal potential. Similarly, $g_L$ and $V_l$ are the effective maximal conductance and reversal potential of any other currents, which are assumed to be independent of the opening and closing of ion channels, and $I$ is an external curent. The four quantities $(g_{\rm Na},g_L,V_{\rm Na},V_L)$ are taken to be constants. Since the right-hand side of (\ref{pwv}) is negative for large $x$ and positive for small $x$, it follows that the voltage $x$ is confined to some interval $\Omega =[x_L,x_R]$. The function $F_n(x)$ is clearly continuous and locally Lipschitz. 

In this example the space $\Gamma$ of discrete states is the set of integers $\{n=0,1,\ldots,N\}$ and the Markov chain is given by a birth-death process: 
\begin{equation}
n\underset{\omega_+(n,x)/\epsilon}\rightarrow n+1,\quad n\underset{\omega_-(n,x)/\epsilon}\rightarrow n-1
\end{equation}
with transition rates
\begin{equation}
 \omega_+(x,n)=\alpha(x) (N-n),\quad \omega_-(x,n)=\beta(x)n.
 \end{equation}
The small parameter $\epsilon$ reflects the fact that sodium channels open at a much faster rate than the relaxation dynamics of the voltage \cite{Keener11}.
It follows that the matrix ${\bf A}(x)$ for fixed $x$ is tridiagonal matrix with
\begin{equation}
\fl A_{n-1,n}(x)=\omega_+(x,n-1),\, A_{nn}(x)=-\omega_+(x,n)-\omega_-(n),\, A_{n+1,n}(x)=\omega_-(x,n+1)
\label{Acomp}
\end{equation}
for $n =0,1,\ldots,N$. It is straightforward to show that the Markov chain is ergodic with unique invariant measure (for fixed $n$) given by
\begin{equation}
\label{ss1}
\rho(x,n)=\frac{N!}{(N-n)!n!}a(x)^nb(x)^{N-n},
\end{equation}
with
\begin{equation}
a(x)=\frac{\alpha(x)}{\alpha(x)+\beta(x)}, \quad b(x)=\frac{\beta(x)}{\alpha(x)+\beta(x)}.
\end{equation}

The above stochastic hybrid system satisfies all of the conditions specified in section 2. Hence, the law of large numbers implies that
in the mean-field limit $\epsilon \rightarrow 0$, we obtain the deterministic kinetic equation
\begin{equation}
\label{bbbrate}
\frac{dx}{dt}=\overline{F}(x)\equiv a(x)f(x)-g(x)
\end{equation}
where
\[
a(x)=\langle n\rangle/N,\quad \langle n\rangle = \sum_{n=1}^{N}n\rho(x,n),
\]
and $\rho$ is the stationary density (\ref{ss1}). One of the features of the averaged model is that it can exhibit bistability for a range of physiologically reasonable parameter values. This is illustrated in Fig. 1, where we plot the deterministic potential $U(x)=-d\overline{F}/dx$ as a function of $x$. Here $x_-$ represents a resting state of the neuron, whereas $x_+$ represents an active state; noise-induced transitions from $x_-$ to $x_+$ can be interpreted in terms of the initiation of a spontaneous action potential. Elsewhere, WKB methods and matched asyptotics have been used to calculate the MFPT to escape from $x_-$. Here, we will focus on the quasipotential and its relation to the Perron eigenvalue 

%%%%%%%%%%%%%%%%%%%%%%%%%%%%%%%%%%%%%%%%%%%%%%%%%%%%%%%%%%%%%%%%%%%%%%%%%%%%%%%%
%%%%%%%%%%%%%%%%%%%%%%%%%%%%%%%%%%%%%%%%%%%%%%%%%%%%%%%%%%%%%%%%%%%%%%%%%%%%%%%%

%%%%%%%%%%%%%%%%%%%%%%%%%%%%%%%%%%%%%%%%%%%%%%%%%%%%%%%%%%%%%%%%%%%%%%%%%%%%%%%%
\begin{figure}[t!]
\begin{center}
\includegraphics[width=12cm]{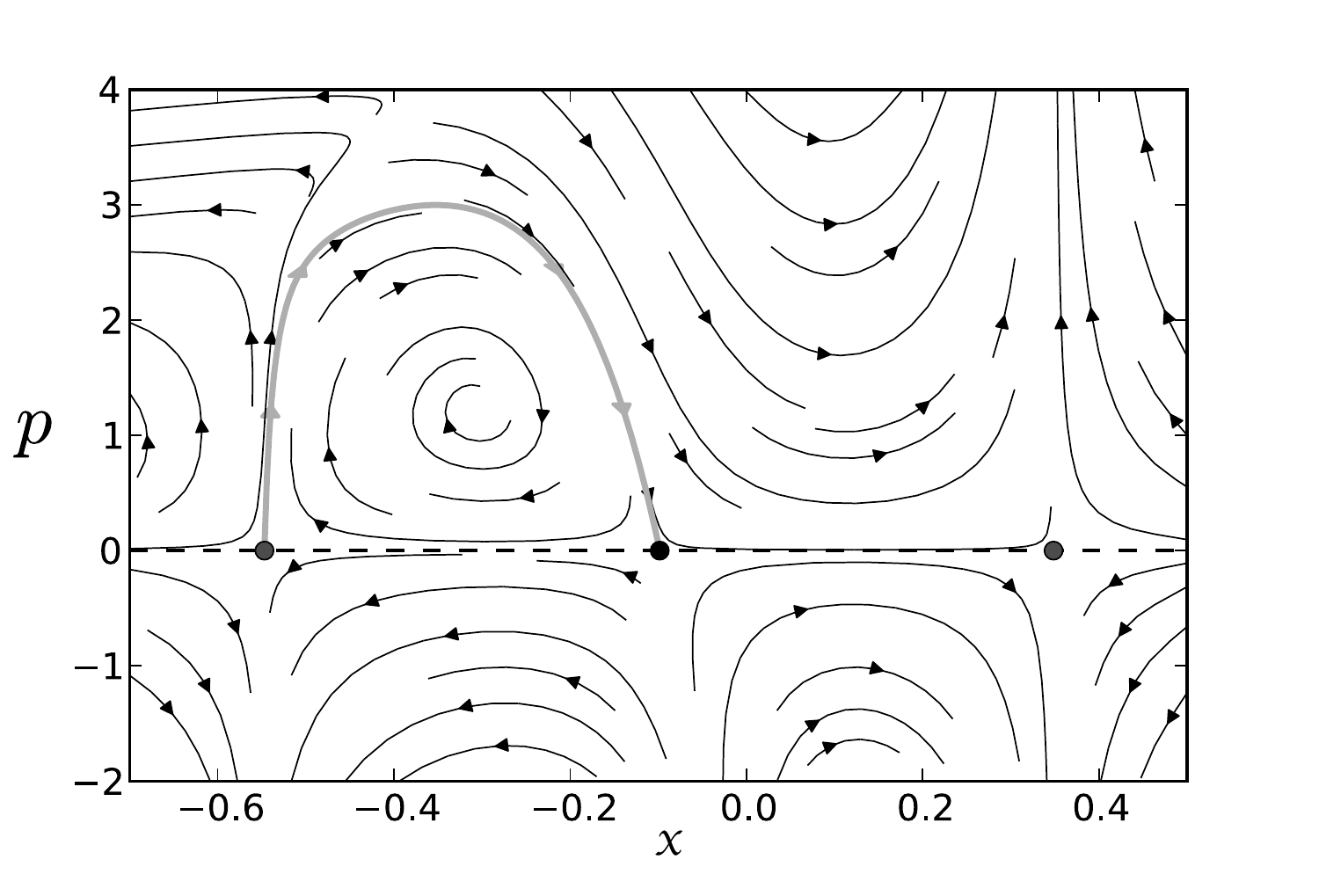}
\caption{\small Phase portrait of Hamilton's equations of motion for the ion channel model with Hamiltonian given by the Perron eigenvalue (\ref{ll}). The zero energy solution representing the maximum likelihood path of escape is shown as the gray curve. Dimensionless units are used with same parameter values as Fig. \ref{wellx} and $I=0$.}
\label{px}
\end{center}
\end{figure}
%%%%%%%%%%%%%%%%%%%%%%%%%%%%%%%%%%%%%%%%%%%%%%%%%%%%%%%%%%%%%%%%%%%%%%%%%%%%%%%%

Substituting the explicit expressions for ${\bf A}$ and $F_n(x)$ into the adjoint equation (\ref{Permz}), yields the following  equation for the Perron eigenvalue $\lambda$ and the left eigenvector ${\bf z}$:
\begin{eqnarray}
&&(N-n+1)\alpha(x) z_{n-1} -[\lambda+n\beta(x)+(N-n)\alpha(x) ]z_n\nonumber
\\ && \quad +(n+1)\beta(x) z_{n+1} =-p \left (\frac{n}{N}f(x)-g(x)\right )z_n
\end{eqnarray}
Consider the trial solution \cite{Bressloff14a}
\begin{equation}
z_n(x,p)=\frac{\Gamma(x,p)^n}{(N-n)!n!},
\end{equation}
which yields the following equation relating $\Gamma$ and $\lambda$:
\begin{eqnarray*}
&\frac{n\alpha}{\Gamma}+\Gamma \beta(N-n)-\lambda    -n\beta-(N-n)\alpha =-p\left (\frac{n}{N}f-g\right ).
\end{eqnarray*}
Collecting terms independent of $n$ and terms linear in $n$ yields the pair of equations
\begin{equation}
\label{pp}
p=-\frac{N}{f(x)}\ \left (\frac{1}{\Gamma(x,p)}+1 \right )\left (\alpha(x)-\beta(x) \Gamma(x,p) \right ),
\end{equation}
and
\begin{equation}
\label{ll}
\lambda(x,p)=-N(\alpha(x)-\Gamma(x,p) \beta(x))-pg(x).
\end{equation}
Eliminating $\Gamma$ from these equation gives
\begin{eqnarray*}
p&=\frac{1}{f(x)}\left ( \frac{N\beta(x)}{\lambda(x,p)+N\alpha(x)+pg(x)}+1\right ) (\lambda(x,p)+pg(x))
\end{eqnarray*}
This yields a quadratic equation for $\lambda$ of the form
\begin{equation}
\label{ll0}
\lambda^2+\sigma(x)\lambda-h(x,p)=0.
\end{equation}
with
\begin{eqnarray*}
\sigma(x)&=(2g(x)-f(x))+N(\alpha(x)+\beta(x)),\\ h(x,p)&=p[-N\beta(x) g(x)+(N\alpha(x)+pg(x))(f(x)-g(x))].
\end{eqnarray*}
Given the Perron eigenvalue, we can determine the quasipotential $\Phi(x)$ by solving the
 Hamiton-Jacobi equation $\lambda(x,\partial_x\Phi)=0$. Equation (\ref{ll0}) then yields the reduced Hamilton-Jacobi equation 
\begin{equation}
h(x,\partial_x\Phi_0)=0.
\end{equation}
The latter is precisely the equation for the quasipotential previously derived using WKB methods \cite{NBK13}. It has the following pair of solutions for $\Phi_0'=\partial_x\Phi_0$:
\begin{equation}
\Phi_0'=0 \mbox{ and } \Phi_0'(x)=- N\frac{\alpha(x) f(x)-(\alpha(x)+\beta)g(x)}{g(x)(f(x)-g(x))}.
\end{equation}
The trivial solution $\Phi_0=$ constant occurs along deterministic trajectories, which converge to the fixed point, whereas the non-trivial solution for $\Phi_0(x)$ occurs along the most likely escape trajectories. In Fig. \ref{px} we show solutions to Hamilton's equations in the $(x,p)$-plane, highlighting the zero energy maximum likelihood curve linking $x_-$ and $x_0$. 

\section{Multi-scale stochastic process}

So far we have assumed that the slow process is piecewise deterministic. However, one of the useful features of taking the Lagrangian LDP \cite{fagg08,fagg09} as our starting point is that it is relatively straightforward to extend our analysis to the case where the slow process is a piecewise SDE. First, recall
that the key idea behind the Faggionato et al. Lagrangian construction is that the slow dynamical process coupled to the fast Markov chain on $\Gamma$ rapidly samples the different discrete states of $\Gamma$ according to some non-negative measure $\psi$. In order to extend this construction to a piecewise SDE, it is necessary to take account of the fact that there are now two levels of stochasticity. That is, after averaging the transition rates of the drift and variance of the SDE  with respect to a given measure $\psi$, the resulting system is still stochastic. Since the slow system operates in a weak noise regime, it follows that one can apply an LDP to the slow system for a given $\psi$. The LDP for the full system is then obtained by combining the rate function of the slow system with the supremum rate function for $\psi$. In this final section, we sketch how to extend the analysis to a piecewise SDE. (Approaching large deviation theory for multi-scale stochastic processes in terms of solutions to an eigenvalue problem has also been considered by Feng and Kurtz \cite{Feng}.) As in section \ref{sect:WKB}, we consider only the case $d=1$.

Consider the piecewise Ito SDE
\begin{equation}
\label{Lang}
dX(t)=F_n(X)+\sqrt{\epsilon}\sigma_n(X)dW(t),
\end{equation}
where $n\in \Gamma$ and $W(t)$ is a Wiener process. The drift term $F_n(X)$ and diffusion term $\sigma_n(X)$ are both taken to be Lipschitz.
When the SDE is coupled to the fast discrete process on $\Gamma$, the stochastic dynamics is described by a differential Chapman-Kolmogorov equation. That is, writing
\[\P\{X(t)\in (x,x+dx),\, n(t) =n|x_0,n_0)=\rho(x,n,t) dx,\]
 we have
\begin{equation}
\label{CKmulti}
\fl \frac{\partial \rho(x,n,t)}{\partial t}=-\frac{\partial }{\partial x}[F_n(x)\rho(x,n,t)]+\frac{\epsilon}{2}\frac{\partial ^2}{\partial x^2}[\sigma_n^{2}(x)\rho(x,n,t)]+\frac{1}{\epsilon}\sum_{m}A^{\top}_{nm}(x)\rho(x,m,t).
\end{equation}
We define the measure space $\meas_+([0,T])$ as before, but now modify the definition of the subspace  ${\mathcal Y}_{x_0}\subset C([0,T]) \times \meas_+([0,T])^\Gamma$ by taking it to be the set of stochastic trajectories satisfying
\begin{eqnarray}
dX(t)&=\sum_{n=1}^K\psi_n F_n(X)+\sqrt{\epsilon}\sqrt{\sum_{n=1}^K \psi_n(t)\sigma_n^2(X)}\, dW(t)\nonumber \\
& \equiv F(X,\psi)+\sqrt{\epsilon}\sigma^2(X,\psi)dW(t) 
\label{SDEj}
\end{eqnarray}
for $\psi\in \meas_+([0,T])^\Gamma$. Such a space contains the set of trajectories of the SDE (\ref{Lang}) with $\psi_n(t)$ given by equation (\ref{pip}) and $n(t)$ evolving according to the Markov chain on $\Gamma$. Consider a particular realization of the Wiener process $W(t)$ on $[0,T]$, which is independent of $\{x(t),\psi(t)\}_{[0,T]}$. For a given realization, one can write down an LDP along identical lines to the case of a piecewise deterministic system, see equation (\ref{LDP2})). Assuming that we can then average with respect to the Wiener process, we obtain an additional contribution to the rate function so that
\begin{eqnarray}
\label{eq:rateM2}
 &J_T(\{x(t)\}_{t\in [0,T]}) =\underset{\{\psi(t)\}_{t\in [0,T]} }\inf  J_T(\{(x(t),\psi(t))\}_{t\in [0,T]}).
\end{eqnarray}
with 
\begin{equation}
\label{eq:rateM}
\fl J_{T}(\{(x(t),\psi(t))\}_{t\in [0,T]})=\int_0^T j(x(t),\psi(t))\,dt+A_T(\{(x(t),\psi(t))\})_{t\in [0,T]})
\end{equation}
for $j(x,\psi)$ given by equation (\ref{eq:j}) and
\begin{equation}
\label{FWact}
A_T(\{x(t),\psi(t)\}_{t\in [0,T]})=\int_0^T\frac{(\dot{x}-F(x,\psi))^2}{2\sigma^2(x,\psi)}dt.
\end{equation}
Equation (\ref{FWact}) is the well known action functional for a one-dimensional SDE in the case of a fixed, time-independent $\psi$ \cite{Freidlin98}. 
One can evaluate the rate function (\ref{eq:rateM2}) along similar lines to the case of a piecewise deterministic system (see appendix A.3) to obtain a classical action with Hamiltonian given by the Perron eigenvalue $\lambda(x,p)$ of the linear equation
\begin{equation}
\fl \sum_{m}A_{nm}(x)R_m(x,p)+(pF_n(x)+p^2\sigma_n^2/2)R_n(x,p)=\lambda(x,p) R_n(x,p).
\end{equation}

\section*
{Acknowledgements}
We thank Jay Newby for drawing our attention to the work of Yuri Kifer \cite{Kifer92,Kifer09} and producing Fig. 2. The work was conducted while PCB was visiting the NeuroMathComp group where he holds an INRIA International Chair. PCB was also partially supported by the NSF (DMS-1120327). OF was partially supported by the European Union Seventh
Framework Programme (FP7/2007-2013) under grant agreement no. 269921 (BrainScaleS), no. 318723
(Mathematics), and by the ERC advanced grant NerVi no. 227747.
\setcounter{equation}{0}
\renewcommand{\theequation}{A.\arabic{equation}}

\appendix
\section{Classical action from the Large Deviation Principle}
\label{app:A}
In this appendix we present the details of deriving the classical action from the LDP given by equation (\ref{LDP2}). 
We first summarize a few useful properties of $j(x,\psi)$ defined by equation (\ref{eq:j}). Since the double sum in equation (\ref{eq:j}) excludes diagonal terms, we introduce the set $\Gamma_\Delta=\Gamma \times \Gamma \backslash \Delta$ where $\Delta$ is the diagonal of $\Gamma \times \Gamma$. Let $c_{nm}=\psi_nW_{nm}(x)$. Suppose that $\psi$ is a strictly positive measure, $\psi_n>0$ for all $n\in \Gamma$. It then follows from the properties of the transition matrix ${\bf W}$ that the mapping $c: \Gamma_\Delta \to [0,\infty)$  is {\em irreducible} in the sense that, for all  $n \neq m \in \Gamma$, there exists a finite sequence $n_1,\,n_2,\cdots,n_k$ such that $n_1=n$, $n_k=m$ and $c_{n_in_{i+1}} > 0$ for $i=1,\cdots,k-1$. Define the mapping $\mI: [0,\infty)^{\Gamma_\Delta} \to \R$ as
\begin{equation}
\label{sym0}
\mI(c)=\sup_{\mR \in (0,\infty)^\Gamma} \hat{\mI}(c,\mR),\quad \hat{\mI}(c,\mR)=\sum_{(n,m) \in \Gamma_\Delta} c_{nm}(1-R_{m}/R_n)
\end{equation}
The following lemma is proven in \cite{fagg08,Kifer09} and follows in part from material found in \cite{rockafellar:70}.
\paragraph{Lemma 1:}
The function $\mI$ is convex and continuous and takes its values in $[0,\infty)$. Moreover, for each irreducible $c $, the supremum on $[0,\infty)^\Gamma$ of the function $\hat{\mI}(c,\cdot)$ is a maximum and is the unique solution of the set of equations
\begin{equation}
\label{sym}
\sum_{m \in \Gamma} c_{nm} \frac{R_{m}}{R_n}=\sum_{m \in \Gamma} c_{mn} \frac{R_{n}}{R_{m}} \quad n \in \Gamma
\end{equation}
under the normalization $\sum_{n \in \Gamma}{R}_n=1$.

\subsection{Evaluating the supremum} We proceed by introducing the following Ansatz regarding the solution $\mR=(R_n,n\in \Gamma)$ of the variational problem (\ref{eq:j}) for fixed $x$ and strictly positive measure $\psi$, namely, that it is an eigenvector of the following matrix equation:
\begin{equation}
\label{Per0}
{\bf A}(x)\mR+\q \circ \mR =\lambda \mR
\end{equation}
for some bounded vector $\q=(q_n,n\in\Gamma)$. We have used the notation introduced in (\ref{eq:ab}).
%Here, for any ${\bf a},{\bf b}\in \R^K$,
% \[[ {\bf a}\circ{\bf b}]_n\equiv [{\rm diag}({\bf a}){\bf b}]_n=a_nb_n.\] 

Note that we are free to shift the vector $\q$ by a constant since, under the transformation $q_n\rightarrow q_n-c$, the eigenvalue shifts by $\lambda\rightarrow \lambda-c$ and the eigenvector is unchanged. That is, for fixed $x$,
\begin{equation}
\lambda(x,\q-c{\bf 1}_K)=\lambda(x,\q)-c.
\end{equation}
For simplicity, we choose $c=q_K$ and take ${\bf Q}=(q_1,q_2,\ldots,q_{K-1},0)$ so that $\lambda(x,{\bf Q})=\lambda(x,\q)-q_K$ and $\mR=\mR(x,{\bf Q})$ are solutions of the matrix equation
\begin{equation}
\label{Per}
{\bf A}(x)\mR+{\bf Q} \circ \mR =\lambda \mR
\end{equation}
There are $K-1$ independent variables, $q_n$, $n=1,\ldots,K-1$.
One of the crucial features of the above Ansatz is that we can then ensure $R_n=R_n(x,{\bf Q})$ is strictly positive by using the Perron-Frobenius theorem and taking $\lambda=\lambda(x,{\bf Q})$ to be the Perron eigenvalue. Indeed, choosing $\kappa$ such that
\begin{equation}\label{kap}
\kappa>\underset{n=1,\ldots K-1}\max\{|q_n|\},
\end{equation}
it is clear that the new matrix $\mA+{\rm diag}(\mathbf{Q})+\kappa \mathbf{I}_K$ is irreducible and positive. According to the Perron-Frobenius theorem, it has a unique strictly positive eigenvector $\z$ with the usual normalization $\sum_m R_m=1$, and this eigenvector is also an eigenvector of $\mA+{\rm diag}(\mathbf{Q})$ (with a shifted eigenvalue, though).
%Proceeding along similar lines to equation (\ref{Wstar}), we consider the modified eigenvalue equation
%\begin{equation}
%\label{PerW}
%\widehat{\bf A}^{\top}(x)\z+{\bf Q} \circ \z =\widehat{\lambda} \z
%\end{equation}
%with $\widehat{\mathbf{A}}(x)$ defined by equation (\ref{Wstar}) and $\widehat{\lambda}(x,{\bf Q})=\lambda(x,{\bf Q})+W^*$.
%We  can apply the Perron-Frobenius Theorem directly to the matrix operator on the left-hand side of equation (\ref{PerW}) provided that we choose $\kappa$ such that
%\begin{equation}
%\label{kap}
%\kappa>\underset{n=1,\ldots K-1}\max\{|q_n|\}.\end{equation}
%This then allows us to extend the Perron-Frobenius Theorem to equation (\ref{Per}) by shifting the Perron eigenvalue.

Applying equation (\ref{Per}) to the double sum in equation (\ref{eq:j}), which holds for strictly positive $\psi$, we find
\begin{eqnarray*}
\fl  \sum_{m,n=1}^K\psi_n W_{nm}(x)\left [1-\frac{R_m}{R_n}\right ]&=\sum_{n=1}^K\psi_n \left [ \sum_{m=1}^K W_{nm}(x)-\frac{1}{R_n}\sum_m (A_{nm}(x)+\delta_{nm})R_m \right ]\\
 &=\sum_{n=1}^K \psi_n[-\lambda+Q_n]\\ &
 =\sum_{n=1}^{K-1} q_n \psi_n -\lambda(x,{\bf Q}),
 \end{eqnarray*}
 because of (\ref{eq:AR}), (\ref{Per}), and since $\sum_{m=1}^K\psi_m=1$ and $Q_K=0$.
In order to ensure that we have found the true supremum with respect to ${\bf z}$, we use equation (\ref{Per}) to show that the supremum satisfies equation (\ref{sym}). Its left-hand side is
\begin{eqnarray*}
\fl  \sum_{m=1}^K\psi_n W_{nm}(x)\frac{R_m}{R_n}&= \frac{\psi_n}{R_n} \sum_{m=1}^KW_{nm}(x)R_m  = \frac{\psi_n}{R_n}(\lambda R_n -Q_n R_n)+R_n=\psi_n (\lambda-Q_n+1).
 \end{eqnarray*}
 The right-hand side reads
 \[
 R_n \sum_{m=1}^K \frac{\psi_m}{R_m} W_{mn}(x)=R_n \sum_{m=1}^K \frac{\psi_m}{R_m} (A_{mn}(x)+\delta_{mn})
 \]
Dividing the left and right-hand sides  by $R_n$,  we deduce that 
\begin{equation}
\label{Rz}
\psi_n=R_n(x,{\bf Q})z_n(x,{\bf Q})
\end{equation}
for each $n=1,\ldots,K$, where $\z =(z_n,n\in \Gamma)$ is the corresponding unique strictly positive eigenvector of the adjoint linear equation
\begin{equation}
  {\bf A}(x)^\top \z+{\bf Q} \circ \z =\lambda \z.
  \label{AR}
  \end{equation}
such that  $\langle \z,\, \mR \rangle=\z^{\top}{\bf R}=1$. Equation (\ref{Rz}) ensures that $\psi$ is a strictly positive measure and that $\sum_{m=1}^K\psi_m=1$. Assuming that such a solution exists, the corresponding function $j(x,\psi)$ is given by
\begin{equation}
\label{jar}
j(x,\psi)=\sum_{n=1}^{K-1} q_n \psi_n-\lambda(x,{\bf Q}),
\end{equation}
where we have dropped the constant term.

In the given variational problem there are $K-1$ independent variables $\psi_n,\, n=1,\ldots,K-1$ with $\psi_K=1-\sum_{n=1}^{K-1}\psi_n$. Similarly, there are $K-1$ independent variables $q_n,n=1,\ldots, K-1$. Therefore, equation (\ref{Rz}) determines a mapping between the sets $\{q_n,\, n=1,\ldots, K-1\}$ and $\{\psi_n,\, n=1,\ldots, K-1\}$. It remains to show that there exists a unique solution ${\bf q}$ for each $\psi \in \meas_+([0,T])^\Gamma$, that is, the mapping is invertible. Differentiating equation (\ref{AR}) with respect to $q_m$, $m=1,\ldots,K-1$, yields the inhomogeneous linear equation
 \begin{eqnarray}
\label{PerL}
\fl {\bf L}(x,{\bf Q})\frac{\partial \z}{\partial q_m}\equiv \left [{\bf A}^\top(x)+{\rm diag}({\bf Q}) -\lambda {\bf I}_K\right ]\frac{\partial \z}{\partial q_m}=\frac{\partial \lambda}{\partial q_m}{\z}-z_m{\bf e}_m,
\end{eqnarray}
where $({\bf e}_m)_n=\delta_{mn}$. Multiplying both sides of (\ref{PerL}) on the left with $\mR^\top$ and using (\ref{Per}) we obtain
%Since the matrix operator ${\bf L}(x,{\bf Q})$ has a one-dimensional null space spanned by the vector ${\bf R}(x,{\bf Q})$ with adjoint vector $\z(x,{\bf Q})$,  the Fredholm alternative implies that
\begin{equation}
\label{lzR}
\frac{\partial \lambda(x,{\bf Q})}{\partial q_m}=R_m(x,{\bf Q})z_m(x,{\bf Q}),\quad m=1,\ldots,K-1.
\end{equation}
Since ${\bf R}$ and ${\bf z}$ are strictly positive, $\lambda(x,{\bf Q}) $ is a monotonically increasing function of the $q_m$.
Moreover, equations (\ref{Per0}) and (\ref{AR}) imply that, in the limit $q_l\rightarrow \infty$ with all other $q_m$ fixed, $R_l,z_l\rightarrow 1$ and $\partial\lambda/\partial q_l \rightarrow 1$. On the other hand, if $q_l\rightarrow -\infty$ then $R_l,z_l\rightarrow 0$ and the Perron eigenvalue becomes independent of $q_l$ with $\partial\lambda/\partial q_l \rightarrow 0$. Hence, by continuity, for each $\psi\in \meas_+([0,T])^\Gamma$ there exists a vector ${\bf Q}$ such that $\psi_n=R_n(x,{\bf Q})z_n(x,{\bf Q})$ for all $n=1,\ldots,K-1$. For such a solution to be unique, the inverse function theorem implies that the Jacobian of the transformation must be invertible. Differentiating equation (\ref{lzR}) with respect to $q_n$ shows that the Jacobian is equivalent to the Hessian of $\lambda$ with respect ${\bf Q}$, since
\begin{equation}
\label{bo}
D_{mn}(x,{\bf Q})\equiv\frac{\partial R_m(x,{\bf Q})z_m(x,{\bf Q})}{\partial q_n}=\frac{\partial^2 \lambda(x,{\bf Q})}{\partial q_m\partial q_n}
\end{equation}
for all $m,n=1,\ldots,K-1$. This also establishes that the Jacobian is a symmetric matrix with real eigenvalues. Invertibility follows from the convexity of the function ${\mathcal J}(c)$ defined by equation (\ref{sym0}). That is, differentiating 
\[j(x,\psi)= {\mathcal J}(c), \quad c_{nm}=\psi_nW_{nm}(x),\]
with respect to $\psi$ for fixed $x$ gives
\begin{equation}
\label{ys}
W_{nm}(x)W_{mm'}(x)\frac{\partial^2 {\mathcal J}(c)}{\partial c_{nm}\partial c_{mm'}}=\frac{\partial^2 j(x,\psi)}{\partial\psi_n\partial \psi_m}
\end{equation}
with $n\neq m'$ and $m\neq m'$. On the other hand, differentiating equation (\ref{jar}) for $j(x,\psi)$ with respect to $\psi_n$ gives
\[
\frac{\partial j(x,\psi)}{\partial\psi_n}=q_n+\sum_{j=1}^{K-1}  \left [\psi_j-\frac{\partial \lambda}{\partial q_j}\right ]\frac{\partial q_j}{\partial \psi_n}=q_n,
\]
because of (\ref{Rz}) and (\ref{lzR}). And so
\[
\frac{\partial^2 j(x,\psi)}{\partial \psi_n\partial \psi_m}=\frac{\partial q_n}{\partial \psi_m} =[{\bf D}^{-1}]_{nm}.
\]
Hence, convexity of ${\mathcal J}(c)$ together with irreducibility of the non-negative transition matrix ${\bf W}$ means that the Jacobian is invertible and positive definite.

In summary, we have shown that for a strictly positive measure $\psi$ and fixed $x$, a unique solution $q_n=q_n(x,\psi)$ exists for all $n=1,\ldots,K-1$ and we have solved the first variational problem by identifying $\z$ with the unique (up to scalar multiplication), strictly positive eigenfunction of the matrix equation (\ref{Per}). Now suppose $\psi$ is a non-negative rather than a strictly positive measure measure, that is, $\psi_m=0$ for at least one state $m\in \Gamma$. In this case $c_{nm}=\psi_nW_{nm}(x)$, $n\neq m$, is not irreducible and lemma 1 no longer applies. However, as proven by Faggionato et al. \cite{fagg08}, the function ${\mathcal J}(c)$ is continuous with respect to $c$. Hence, assuming that the form of the rate function (\ref{LDP2}) still holds (which isn't necessarily true), we can take a sequence of strictly positive measures $\psi^{(l)}$ on $\Gamma$ such that $\psi_n^{(l)}\rightarrow \psi_n$ for each $n\in \Gamma$. This implies that (for fixed $x$)
\[ c[\psi^{(l)}]\rightarrow c[\psi]\]
and
\[j(x,\psi^{(l)})={\mathcal J}(c[\psi^{(l)}])\rightarrow {\mathcal J}(c[\psi])=j(x,\psi).\]
so that one can extend equation (\ref{jar}) to non-negative measures by taking
\[q_n(x,\psi)=\lim_{l\rightarrow \infty} q_n(x,\psi^{(l)}).\]

\paragraph{Example.}
 We will illustrate the Perron eigenvalue solution to the supremum variational problem by considering an example for $K=2$. Let us take the transition matrix to be
 \[
{\bf W}= \left (\begin{array}{cc} 1/2 & 1/3 \\ 1/2 & 2/3 \end{array}
 \right )
 \]
 Consider the eigenvalue equation
 \[{\bf W}{\bf R} +\q \circ{\bf R}=\lambda {\bf R},\]
 where we have absorbed the diagonal terms $\sum_{k=1,2}W_{km}$ into the definition of $q_m$.
 The resulting characteristic equation is a quadratic in $\lambda$ and the leading or Perron eigenvalue is given by
 \begin{eqnarray*}
 \lambda=\frac{q_1+q_2}{2}+\frac{7}{12}+\frac{1}{2}\sqrt{(q_1-q_2)^2-(q_1-q_2)/3+25/36}.
 \end{eqnarray*}
 It follows that
 \begin{eqnarray*}
 \psi_1&\equiv \frac{\partial \lambda}{\partial q_1}=\frac{1}{2}+f(q_1-q_2)\\
 \psi_2&\equiv \frac{\partial \lambda}{\partial q_2}=\frac{1}{2}-f(q_1-q_2)
 \end{eqnarray*}
 with
 \[f(q)=\frac{1}{4}\frac{2q-1/3}{\sqrt{q^2-q/3+25/36}}.\]
 Note that $\psi_1+\psi_2=1$ as required.
 The function $f(q)$ is a monotonically increasing function of $q$ with $f(-\infty)=-1/2$ and $f(\infty)=1/2$. Thus, one can find a unique, finite value of $q=q_1-q_2$ for all $\psi_1\in (0,1)$, that is, for all strictly positive $\psi$. In the case of a non-negative $\psi$ with $\psi_1=0$ or $\psi_2=0$, we have $q\rightarrow \pm \infty$.

\subsection{Evaluating the infimum}
The next step is to substitute for $j(x,\psi)$ in the rate function (\ref{eq:rate2}), which gives
\begin{equation}
\fl J_T(\{x(t)\}_{t\in [0,T]})=\underset{\psi: \dot{x}=\sum_{n=1}^K\psi_n \mF_n(x)} \inf  \int_0^T\left [\sum_{n=1}^{K-1}q_n(t)\psi_n(t)-\lambda(x(t),{\bf Q}(t))\right ]dt,
\end{equation}
with $q_n(t)=q_n(x(t),\psi(t))$ and $\sum_{m=1}^K\psi_m =1$. In order to solve this variational problem, we
introduce a $d$-dimensional Lagrange multiplier ${\bm \mu}(t)$ and set
\begin{equation}
{J}_T(\{x(t)\}_{t\in [0,T]}) =\underset{\psi,{\bm \mu}} \inf\,  {S}[x,{\bm \mu},\psi],
\end{equation}
where
\begin{eqnarray}
{S}[x,{\bm \mu},\Psi]&= \int_0^T \left [\sum_{n=1}^{K-1}q_n(t)\psi_n(t)-\lambda(x(t),{\bf Q}(t))\right . \\
&\qquad \qquad \left . +\langle {\bm \mu}(t),\,\dot{x}-\sum_{n=1}^{K-1}[\mF_n(x)-\mF_K(x)]\psi_n(t)-\mF_K(x)\rangle\right ]dt,\nonumber
\end{eqnarray}
and we have imposed the constraint $\sum_{m=1}^K\psi_m =1$.
The variational problem can now be expressed in terms of functional derivatives of $S$:
\begin{equation}
\frac{\delta {S}}{\delta {\bm \mu}(s)}=\dot{x}(s)-\sum_{n=1}^K\mF_n(x(s))\psi_n(s)=0,
\end{equation}
and
\begin{equation}
\fl \frac{\delta {S}}{\delta \psi_m(s)}=\sum_n\frac{\partial q_n}{\partial \psi_m}\psi_n(s)+q_m(s)-\sum_{n=1}^{K-1}\frac{\partial \lambda}{\partial q_n}\frac{\partial q_n}{\partial \psi_m} -\langle {\bm \mu}(s),\,\mF_m(x(s))-\mF_K(x(s)) \rangle=0
\end{equation}
for $m=1,\ldots K-1$.
Combining with equations (\ref{Rz}) and (\ref{lzR}), we obtain the following solution to the variational problem in terms of $\bm \mu$:
\begin{eqnarray}
&  q_m=\langle {\bm \mu},\,\mF_m(x)-\mF_K(x) \rangle, \nonumber\\ 
 & \psi_m(x,{\bm \mu})= z_m(x,{\bm \mu})R_m(x,{\bm \mu}) \label{cond}, 
 \end{eqnarray}
for all $m=1,\ldots K-1$, with
 $R_n(x,{\bm \mu}),z_n(x,{\bm \mu})$ the positive eigenvectors of the matrix equations
\begin{equation}
\label{Permu}
 {\bf A}{\bf R}+({\bf F}(x){\bm \mu}) \circ {\bf R}=\lambda {\bf R}
\end{equation}
and its adjoint, respectively. The $K$-dimensional vector ${\mathbf F}(x){\bm \mu}$ is the product of the $K \times d$ matrix ${\mathbf F}(x)$ whose $K$ rows are the $d$-dimensional vectors $\mF_m(x)$, $m=1,\cdots,K$, with the $d$-dimensional vector $\bm \mu$. Here
\[\lambda=\lambda(x,{\bm \mu})\equiv \lambda(x,{\bf Q}|q_m=\langle {\bm \mu},\, \mF_m-\mF_K \rangle,m=1,\ldots,K-1).\]
The corresponding function $j(x,\psi)$ becomes
\[
j(x,\psi)=\sum_{n=1}^K\psi_n(x,{\bm \mu})\langle {\bm \mu},\,\mF_n(x) \rangle-\lambda(x,{\bm \mu}).
\]

The final step is to show that for each $x$, the equation
\[\dot{x}=\sum_{n=1}^{K}\mF_n(x)\psi_n(x,{\bm \mu})\]
is invertible so that the function ${\bm \mu}={\bm \mu}(x,\dot{x})$ exists.
From the inverse function theorem we require that 
\begin{equation}\label{eq:inverse}
\sum_{n=1}^K \frac{\partial \psi_n(x,\mu)}{\partial {\bm \mu}}\mF_n(x) \neq 0,
\end{equation}
for all $x\in \Omega$. Note that the left hand side is a scalar.
Following along identical lines to the analysis of equation (\ref{PerL}), we differentiate the linear equation (\ref{Permu}) with respect to ${\bm \mu}$ to give
\begin{equation}
\label{goog}
{\bf L}(x,{\bm \mu})\frac{\partial {\bf R}}{\partial {\bm \mu}}   ={\bf R}\frac{\partial \lambda}{\partial {\bm \mu}} -\left [{\bf G}_1\circ {\bf R} \cdots {\bf G}_d \circ {\bf R}\right],
\end{equation}
with ${\bf G}_k$, $k=1,\cdots,d$ the $K$-dimensional vector of the $k$th coordinates of the $K$ vectors $\mathbf{F}_m(x)$, $m=1,\cdots,K$, and ${\bf L}(x,{\bm \mu})=\left . {\bf L}(x,\p)\right |_{p_m=\langle {\bm \mu},\, \mF_m \rangle}$. Using the same arguments as previously we obtain
%Since the matrix ${\bf L}$ has a one-dimensional null-space spanned by the adjoint eigenvector ${\bf R}(x,{\bm \mu})$, the Fredholm alternative together with the normalization $\z^{\top}{\bf R}=1$ implies that
\begin{equation}
\label{doty}
\frac{\partial \lambda(x,{\bm \mu})}{\partial {\bm \mu}}^\top=\sum_{n=1}^K\mF_n(x)z_n(x,{\bm \mu})R_n(x,{\bm \mu})
\end{equation}
Combining (\ref{bo}), (\ref{cond}), (\ref{doty}) and (\ref{eq:inverse}) we require
\begin{eqnarray*}
\fl \frac{\partial^2 \lambda(x,{\bm \mu})}{\partial {\bm \mu}^2}&\equiv \sum_{m,n=1}^{K-1}\left .\frac {\partial^2 \lambda(x,{\bf Q})}{\partial q_m\partial q_n}\right |_{q_m=\langle {\bm \mu},\, \mF_m(x)-\mF_K(x) \rangle} \langle \mF_m(x)-\mF_K(x),\, \mF_n(x)-\mF_K(x) \rangle  \\
\fl &\neq 0. 
\end{eqnarray*}
If we define $X$ to be the $(K-1) \times d$ matrix whose rows are the $d$-dimensional vectors $\mF_m(x)-\mF_K(x)$, $m=1,\cdots,K-1$ we have
\[
\fl \frac{\partial^2 \lambda(x,{\bm \mu})}{\partial {\bm \mu}^2}={\rm Trace}(\mathbf{D}XX^\top)={\rm Trace}(X^\top\mathbf{D}X)=\sum_{k=1}^d X_k^\top \mathbf{D} X_k
\]
The second equality is true because of the properties of the Trace operator. The $(K-1)$-dimensional vectors $X_k$, $k=1,\cdots,d$ are the column vectors of the matrix $X$. All terms in the right hand side are positive since the Jacobian $\mathbf{D}$ has been shown to be positive definite. At least one of the vectors $X_k$ is non zero since $\mF_m(x)\neq \mF_K(x)$ for at least one $m\neq K$ and the corresponding term in the sum in the right hand side of the last equality is strictly positive.

Finally, from equation (\ref{kap}), we require ${\bm \mu}$ to be bounded, that is, there exists a $\kappa$ for which
\[\kappa > \underset{m=1,\ldots,K}\max\{|\langle {\bm \mu},\, \mF_m(x) \rangle |\}\]
for all $x\in \Omega$.

\subsection{Extension to a multi-scale process}
The above analysis can be extended to the LDP (\ref{eq:rateM2}) for a multi-scale stochastic process. In particular, using equation (\ref{jar}) we have  
\begin{equation}
 J_T(\{x(t)\}_{t\in [0,T]})=\underset{\{\psi(t)\}_{t\in [0,T]} }\inf S[x,\psi]
\end{equation}
where
\begin{eqnarray}
\fl S[x,\psi]=\int_0^T\left [\sum_{n=1}^{K-1}q_n(t)\psi_n(t)-\lambda(x(t),{\bf Q}(t)) +\frac{(\dot{x}-\sum_{n=1}^{K}\psi_nF_n(x))^2}{2\sum_{n=1}^{K-1}\psi_n \sigma_n^2} \right ]dt
\end{eqnarray}
with $q_n(t)=q_n(x(t),\psi(t))$, $\sum_{m=1}^K\psi_m =1$ and $\lambda$ the Perron eigenvalue of equation (\ref{Per0}). Taking the infimum with respect to $\psi_k$, $k=1,\ldots, K-1$, gives
\begin{eqnarray}
\fl &0=\frac{\delta S}{\delta \psi_k}=q_k-\frac{(\dot{x}-\sum_{n=1}^K\psi_nF_n(x))[F_k(x)-F_K(x)]}{\sum_{n=1}^K\psi_n \sigma_n^2}  \nonumber \\
\fl &\qquad -\frac{(\dot{x}-\sum_{n=1}^K\psi_nF_n(x))^2[\sigma_k^2(x)-\sigma_K^2(x)]}{2[\sum_{n=1}^K\psi_n \sigma_n^2]^2} +\sum_{n=1}^{K-1}\left (\psi_n -\frac{\partial \lambda}{ \partial q_n}\right )\frac{\partial q_n}{\partial \psi_k}
\end{eqnarray}
Introducing the new variables
\begin{equation}
{\bm \mu}=\sum_{n=1}^K\psi_nF_n(x),\quad \sigma^2=\sum_{n=1}^K\psi_n\sigma_n^2,\quad p=\frac{\dot{x}-{\bm \mu}}{\sigma^2},
\end{equation}
and noting that $\partial \lambda/\partial q_m=\psi_m$, we have
\begin{equation}
q_k=p[F_k(x)-F_K(x)]+\frac{p^2}{2}[\sigma_k^2(x)-\sigma_K^2(x)]
\end{equation}
and
\begin{eqnarray}
\fl S[x,\psi]&=\int_0^T\left [p\sum_{n=1}^KF_n\psi_n+\frac{p^2}{2}\sum_{n=1}^K\psi_n\sigma_n^2-\lambda(x,p)  +\frac{(\dot{x}-\mu)^2}{2\sigma^2} \right ]dt\\
\fl &=\int_0^T\left [p\mu+\frac{p^2}{2}\sigma^2-\lambda(x,p)  +\frac{p}{2}(\dot{x}-\mu) \right ]dt\\
\fl &=\int_0^T\left [p\dot{x}-\lambda(x,p) \right ]
\end{eqnarray}
with $\lambda$ the Perron eigenvalue for the linear equation (which is independent of $\mu$)
\begin{equation}
\sum_{n}A_{mn}R_n+(pV_m(x)+p^2\sigma_m^2/2)R_m=\lambda R_m.
\end{equation}
Finally $\mu$ and $\sigma^2$ are determined from the identities
\begin{equation}
\fl \mu=\sum_jR_j(x,p)z_j(x,p)V_j(x),\quad \sigma^2 =\sum_jR_j(x,p)z_j(x,p)\sigma^2_j(x).
\end{equation}

\section{Relation with the work of Kifer}

In this appendix, we briefly describe the version of an LDP for stochastic hybrid systems developed in some detail by Kifer \cite{Kifer09}. We first note that Kifer considers a more general dynamical system than presented here, namely, a random evolution on $\Omega \times M \times \Gamma$ with $M$ a compact Riemannian manifold. We will consider 
the more restricted case $M=\emptyset$. Let $(n^\varepsilon_{x_0,n_0}(t),x^\varepsilon_{x_0,n_0}(t))$ denote a particular solution of the stochastic hybrid system evolving according to equation (\ref{hs}) with initial conditions $(x_0,n_0)$. Kifer  \cite{Kifer09} proves that the Hamiltonian
\begin{equation}\label{eq:HK}
H(x_0,x,\beta)= \lim_{t \to \infty} \frac{1}{t} \log {\mathbb E} \left[ \exp \beta\int_0^t F_{n^\varepsilon_{x_0,n_0}(s)}(x)\,ds\right]
\end{equation}
exists uniformly in $x_0,\,x \in \bar{\Omega}$ and $|\beta| \leq b$ for some positive constant $b$. Moreover, $H$
is strictly convex in $\beta$ and Lipschitz continuous in the other variables, and does not depend on $n_0$. The Hamiltonian is then used to define the Lagrangian $L(x,x',\alpha)$ as the Fenchel-Legendre transform of $H$,
\[L(x,x',\alpha)=\underset{\beta}\sup\left (\alpha \beta -H(x,x',\beta)\right ),\]
 which is non-negative, (strictly) convex, lower semicontinuous and $H$ is its inverse Fenchel-Legendre transform. Note the quantities that map onto our Hamiltonian and Lagrangian are $H(x,\beta)=H(x,x,\beta)$ and $L(x,\alpha)=L(x,x,\alpha)$.

Consider $C_{0T}$ the space of continuous mappings $[0,T] \to \Omega$. For each $\gamma \in C_{0T}$ note $\gamma_t=\gamma(t)$, $t \in [0,T]$. For each absolutely continuous $\gamma$ its velocity $\dot{\gamma}$ exists for almost all $t \in [0,T]$ and $t \to \dot{\gamma}_t$  is measurable. Hence one can define the action
\[
S_{0T}(\gamma)=\int_0^T L(\gamma_t,\dot{\gamma}_t)\,dt .
\]
It can be shown that $S_{0T}$ is lower semicontinuous with respect to the $sup$ norm on $C_{0T}$, which implies that the set $\Psi^a_{0T}(x)=\{\gamma \in C_{0T}\,:\,\gamma_0=x,\, S_{0T}(\gamma) \leq a\}$ is a closed set. Introducing the the time-rescaled solution $z^\varepsilon_{x_0,n_0}(t)=x^\varepsilon_{x_0,n_0}(t/\varepsilon)$, Kifer establishes the existence of an LDP for the family of random paths $\{z^\varepsilon_{x_0,n_0}\}$: For all $a,\,\delta,\,\lambda >0$ and every $\gamma \in C_{0T}$ such that $\gamma_0=x_0$ there exists $\varepsilon_0=\varepsilon_0(x_0,\gamma,a,\delta,\lambda)$ such that for $\varepsilon < \varepsilon_0$, uniformly in $n_0 \in \Gamma$, we have the lower bound on the open sets
\[
P(\left| z^\varepsilon_{x_0,n_0}-\gamma \right| < \delta) \geq \exp\left\{ -\frac{1}{\varepsilon}(S_{0T}(\gamma)+\lambda) \right\},
\]
and the upper bound on the closed sets
\[
P(\left| z^\varepsilon_{x_0,n_0}-\Psi^a_{0T}(x_0) \right| \geq \delta) \leq \exp\left\{ -\frac{1}{\varepsilon}(a-\lambda) \right\},
\]
Note that the initial condition is taken to be sufficiently far from the boundary $\partial \Omega$ by requiring $\inf_{z\in \partial \Omega}|x_0-z|\geq 2KT$. By the contraction principle one can then obtain an asymptotic estimate of the probability of the rescaled slow process to exit from a given open set $V \subset \Omega$. That is, defining $\tau^\varepsilon_{x_0,n_0}(V)$ to be $\inf\{  t \geq 0: z^\varepsilon_{x_0,n_0}(t) \notin V \}$, then
\[
\fl \lim_{\varepsilon \to 0} P(\tau^\varepsilon_{x_0,n_0}(V) < T)=
\exp \left\{ -\frac{1}{\varepsilon} \inf\{ S_{0t}(\gamma):\, \gamma \in C_{0T},\, t \in [0,T],\, \gamma_0=x_0,\,\gamma_t \notin V \} \right\}
\]

In order to connect with our analysis of the LDP obtained by Faggionato et al \cite{fagg08,fagg09}, we note that the Hamiltonian given by equation (\ref{eq:HK}) can be expressed in terms of a Perron eigenvalue. The first step is to define the transition matrix according to
\[
\fl P(n^\varepsilon_{x_0,n_0}(t+\Delta)=n\,|\,n^\varepsilon_{x_0, n_0}(t)=m,\,x^\varepsilon_{x_0,n_0}(t)=x)=W_{nm}(x) \Delta+o(\Delta)
\]
The infinitesimal generator $\LL^{x_0}$ of the continuous time Markov chain defined by the matrix $\mathbf{W(x_0)}$ is the operator acting on vectors $z \in \R^{|\Gamma|}$ by the formula
\[
(\LL^{x_0}\,z)_k=\sum_{l \in \Gamma} W_{lk}(x_0)(z_l-z_k)
\]
It can then be proven \cite{Kifer92,Kifer09} that $H(x_0,x,\beta)$ is the principal or Perron eigenvalue of the matrix operator $\LL^{x_0}+\beta F_\cdot(x):$ $\R^{|\Gamma|} \to \R^{|\Gamma|}$ acting on vectors $z \in \R^{|\Gamma|}$ by
\begin{equation}
\label{Kper}
\fl ((\LL^{x_0}+\beta F_\cdot(x))z)_k=(\LL^{x_0}\,z)_k+\beta F_k(x) z_k=\sum_{l \in \Gamma} W_{lk}(x_0)(z_l-z_k)+\beta F_k(x) z_k
\end{equation}
The corresponding Lagrangian $L(x_0,x,\alpha)$ is then given in explicit form by
\[
L(x_0,x,\alpha)=\inf_\psi \{ I_{x_0}(\psi)\,:\,\sum_{k \in \Gamma} F_k(x) \psi_k=\alpha \}
\]
where $\psi$ is a probability measure on $\Gamma$, and $I_{x_0}(\psi)$ is defined by
\[
I_{x_0}(\psi)=-\inf_{z > 0} \sum_{k \in \Gamma} \psi_k \frac{(\LL^{x_0} z)_k}{z_k}=\sup_{z>0} \sum_{k,l \in \Gamma} \psi_k W_{lk}(x_0)(1-\frac{z_l}{z_k})
\]
Note that $I_x(\psi)$ is precisely the function $j(x,\psi)$ of equation (\ref{eq:j}), $L(x,\alpha)=L(x,x,\alpha)$ is equivalent to the definition of the Lagrangian appearing in equation (\ref{eq:rate2}), and the eigenvalue equation (\ref{Kper}) for $x_0=x$ is equivalent to equation (\ref{Permz}) with $\mu=\beta$.

\end{document}